\documentclass[12pt]{amsart}
\allowdisplaybreaks
\usepackage{amssymb}
\usepackage{graphicx}
\usepackage{mdwlist}
\usepackage{fullpage}
\usepackage[pdfencoding=auto,colorlinks,linkcolor=,citecolor=]{hyperref}
\usepackage{cite}
\usepackage{dsfont}
\usepackage{eucal}
\usepackage[all,cmtip]{xy}

\makeindex
\numberwithin{equation}{section}

\theoremstyle{plain}
\newtheorem{lemma}[equation]{Lemma}
\newtheorem{theorem}[equation]{Theorem}
\newtheorem{proposition}[equation]{Proposition}
\newtheorem{corollary}[equation]{Corollary}

\theoremstyle{definition}

\newtheorem{definition}[equation]{Definition}
\newtheorem{example}[equation]{Example}
\newtheorem{notation}[equation]{Notation}

\theoremstyle{remark} 
\newtheorem{remark}[equation]{Remark}

\newcommand{\bu}{\bullet}
\newcommand{\Ext}{\mathsf{Ext}}

\renewcommand{\ge}{\geqslant}
\newcommand{\Hom}{\mathsf{Hom}}

\renewcommand{\le}{\leqslant}
\newcommand{\NW}{{\cN_W}}
\newcommand{\kNW}{{\kk\NW}}
\newcommand{\Mat}{\mathsf{Mat}}
\newcommand{\NS}[1]{{\cN_{S_{#1}}}}
\newcommand{\NSn}{\NS{n}}
\newcommand{\kNSn}{{\kk\NSn}}
\newcommand{\op}{{\mathsf{op}}}
\newcommand{\rank}{\mathsf{rank}}
\newcommand{\Tot}{\mathsf{Tot}\,}
\newcommand{\Tr}{\mathsf{Tr}}

\newcommand{\cN}{\mathcal{N}}
\newcommand{\tC}{\tilde{C}}
\newcommand{\td}{\tilde{d}}
\newcommand{\tz}{\tilde z}
\newcommand{\vz}{\check{z}}
\newcommand{\kk}{{\mathbf{k}}}
\newcommand{\bZ}{{\mathbb{Z}}}
\newcommand{\bQ}{{\mathbb{Q}}}

\author{David J. Benson} 
\address{Institute of Mathematics \\ 
University of Aberdeen \\ 
Aberdeen AB24 3UE \\ 
United Kingdom}

\title{The cohomology of the nilCoxeter algebra}

\begin{document}

\begin{abstract}
The nilCoxeter algebra $\NSn$ of the symmetric group $S_n$ is the
algebra over $\bZ$ with generators $Y_i$ ($1\le i\le n-1$),
satisfying the braid relations $Y_iY_{i+1}Y_i=Y_{i+1}Y_iY_{i+1}$,
$Y_iY_j=Y_jY_i$ ($|j-i|\ge 2$), together with the relations $Y_i^2=0$.
We describe an explicit presentation for the cohomology ring 
$Z\cong\Ext^*_{\NSn}(\bZ,\bZ)$, with $n-i$ new generators in degree $i$
for $0< i<n$, and all relations are quadratic. We show that this $\Ext$ ring is $\bZ$-free, 
and that it is a semiprime Noetherian affine polynomial identity (PI)
ring with Poincar\'e series $1/(1-t)^{n-1}$ and PI degree
$2^{n-2}$. For any field of coefficients $\kk$, we show that
$\Ext^*_{\kNSn}(\kk,\kk)$ is $\kk\otimes_\bZ Z$. Similar results hold
for other finite Coxeter types. In the final section we show that $Z$
is a Koszul algebra whose Koszul dual is 
a signed version of the nilcactus algebra, an algebra closely related to the cactus
group.
\end{abstract}



\maketitle

\section{Introduction}

The nilCoxeter algebra 
was studied by
Fomin and Stanley~\cite{Fomin/Stanley:1994a} (1994) as an algebra of divided
difference operators related to the Schubert polynomials, but it
already appears in several earlier papers. See 
Bernstein, Gelfand and Gelfand~\cite{Bernstein/Gelfand/Gelfand:1973a}
(1973) where it is implicit in the proposition on page~5, and
Kostant and Kumar~\cite{Kostant/Kumar:1986a} (1986), where goes under the
name of nil Hecke ring in Section~4.
This algebra
was further investigated by
Khovanov~\cite{Khovanov:2001a}, Bazlov~\cite{Bazlov:2006a}, and Yang~\cite{Yang:2015a} among
others. In particular, it is shown in~\cite{Khovanov:2001a} that the
representation theory of the Weyl algebra
$\bZ\langle x,\partial\rangle/(\partial x - x\partial-1)$ is
categorified by the nilCoxeter algebras of the symmetric groups. 
In this paper, we describe the structure of the cohomology ring of
these algebras. Our main theorem is as follows.

\begin{theorem}\label{th:main}
Let $\NSn$ be the nilCoxeter algebra of the symmetric group
$S_n$ with $n\ge 2$. Then the ring $Z\cong\Ext^*_{\NSn}(\bZ,\bZ)$ is
$\bZ$-free, with Poincar\'e series 
\[ \sum_{i=0}^\infty t^i\rank_\bZ\,Z^i=
  \frac{1}{(1-t)^{n-1}}. \]
The ring $Z$ is given by
generators and relations as follows. There are generators $z_{i,j}$
for $1\le i,j\le n$ of degree $|z_{i,j}|=|j-i|$, satisfying the following
relations.
\begin{enumerate}
\item $z_{j,i}=(-1)^{j-i}z_{i,j}$, $z_{i,i}=1$.
\item If $[i',j']\subseteq [i,j]$ then
$z_{i,j}z_{i',j'} =(-1)^{(j-i)(j'-i')}z_{i+j-i',i+j-j'}z_{i,j}$. 
\item If $[i,j]$ is disjoint from $[i',j']$ then
$z_{i,j}z_{i',j'} =(-1)^{(j-i)(j'-i')}z_{i',j'}z_{i,j}$.
\item If $[i,j]$ and $[i',j']$ overlap, but neither is contained in the other, then
  $z_{i,j}z_{i',j'}=0$.
\end{enumerate}
The ring $Z$ is not commutative, even up to signs, because of relation
{\rm (2)}, but is a semiprime Noetherian 
affine PI ring of PI degree $2^{n-2}$ and Gelfand--Kirillov dimension $n-1$.

For any field $\kk$ of coefficients, we have
$ \Ext^*_{\kNSn}(\kk,\kk)\cong \kk
  \otimes_\bZ Z$, 
a $\kk$-algebra which also has PI degree $2^{n-2}$ and Gelfand--Kirillov dimension $n-1$.
\end{theorem}

The theorem follows from Corollaries~\ref{co:PS} and~\ref{co:Ext} and
Theorem~\ref{th:PI-degree}. 
The method of proof involves producing a free resolution of $\bZ$ as
an $\NSn$-module, that is minimal
when tensored with any field. The resolution comes in the form of the
total complex of an $(n-1)$-fold multiple complex with a free
$\NSn$-module of rank one
at each multidegree, but the differentials are quite complicated to
describe.

Similar results hold for any nilCoxeter algebra.
The generators
correspond to connected oriented subgraphs of finite type 
in the Coxeter diagram. Relation (1) is
for reversing orientation. Relation (2), instead of reflecting the
interval $[i',j']$ in the interval $[i,j]$, conjugates the 
connected subgraph by the
longest word in the Coxeter subgroup for the subdiagram
corresponding to the connected graph. Relation (3) is for disjoint
subgraphs with no joining edge, and relation (4) is for subgraphs that
overlap or have a joining edge without either being
a subgraph of the other. The Poincar\'e
series is equal to $1/(1-t)^r$ where $r$ is the rank of the Coxeter
group. The proofs are
essentially the same, but notationally much more cumbersome to express,
so we limit ourselves to describing the symmetric group case in
detail, and giving the examples of type $B_2$, $G_2$, $I_2(n)$, $B_3$,
$H_3$, and $D_4$ to illustrate what
happens more generally. The same ideas
also work for the super
nilCoxeter algebras of Section~8 of Rosso and
Savage~\cite{Rosso/Savage:2015a};
we describe the case of type super $A_3$ for comparison.

It would be interesting but probably difficult to compute the
structure of the Hochschild cohomology ring of the nilCoxeter
algebras, or at least to determine whether the $\Ext$ ring is finitely
generated as a module over the image of Hochschild cohomology.\bigskip

In the final section, we discuss the following.

\begin{theorem}\label{th:main2}
The ring $Z$ is a graded Koszul algebra whose graded Koszul dual $X$
is again $\bZ$-free, and is a signed version of the nilcactus algebra.
The ring $X$ is given by generators and relations
as follows. There are generators $X_{i,j}$ for $1\le i, j\le n$, of
degree $|X_{i,j}|=|j-i|-1$ satisfying the following relations.
\begin{enumerate}
\item $X_{j,i}=(-1)^{j-i-1}X_{i,j}$, $X_{i,i}=0$.
\item $X_{i,j}^2=0$.
\item If $[i',j']\subseteq [i,j]$ then
  $X_{i,j}X_{i',j'}=(-1)^{(j-i-1)(j'-i'-1)}X_{i+j-i',i+j-j'}X_{i,j}$.
\item If $[i,j]$ is disjoint from $[i',j']$ then
  $X_{i,j}X_{i',j'}=(-1)^{(j-i-1)(j'-i'-1)}X_{i',j'}X_{i,j}$.
\end{enumerate}
\end{theorem}

We give $X$ a triple grading with
$|X_{i,j}|=(-1,|j-i|,\frac{|j-i|(|j-i|+1)}{2})$, 
and describe the spectral sequence of Eilenberg--Moore type 
$X \Rightarrow \NSn$. Although $X$ and $\NSn$ have isomorphic
cohomology rings, no homomorphism between them can induce this
isomorphism because one cohomology ring is formal and the other is not.\bigskip

\noindent
{\bf Convention.} For the sake of the introduction, we have written
all degrees positively. But for Koszul dual algebras, we should really 
index one positively and the other negatively. We choose from now on
to write everything with homological rather than cohomological
grading, so that the degrees in Theorem~\ref{th:main} should be
negated, with $|z_{i,j}|=-|j-i|$, whereas the degrees in
Theorem~\ref{th:main2} should remain positive. Thus differentials
always lower degree.\bigskip

\noindent
{\bf Acknowledgements.} I'd like to thank Pavel Etingof, Kay Jin Lim, and Toby
Stafford for interesting conversations related to this work.

\section{The nilCoxeter algebra}\label{se:Y}

Let $W$ be a finite Coxeter group with Coxeter generators $w_i$ $(i\in I)$, relations
$w_i^2=1$, $(w_iw_j)^{n_{ij}}=1$, length function $\ell$, and longest
element $w_0$.
Note that $w_0^2=1$, so conjugation by $w_0$ is an involutory
automorphism. It permutes the Coxeter generators $w_i$, so it induces
an automorphism of the Coxeter diagram, which is as follows. In types
$B_n$, $C_n$, $D_n$ ($n$ even), $E_7$, $E_8$, $F_4$, $G_2$, 
$H_3$, $H_4$, $I_2(n)$ ($n$ even) it is the identity. In
types $A_n$, $D_n$ ($n$ odd), $E_6$, $I_2(n)$ ($n$ odd)
it is the unique non-trivial involutory automorphism (recall that
$I_2(3)=A_2$, $I_2(4)=B_2$, $I_2(5)=H_2$, $I_2(6)=G_2$). We are particularly interested in the
case $W=S_n$, the symmetric group of degree $n$, with Coxeter type
$A_{n-1}$, but we shall discuss some other types in Section~\ref{se:eg}.

The nilCoxeter algebra $\NW$ over $\bZ$
has generators $Y_i$ $(i\in I)$ with relations $Y_i^2=0$, $Y_iY_j\dots=Y_jY_i\dots$ ($n_{ij}$ terms
on each side). For example, if $W$ is the symmetric group $S_n$  on the
symbols $\{1,\dots,n\}$ then $w_i$ is the permutation $(i,i+1)$ for
$I=\{1,\dots,n-1\}$, and the presentation for $\NSn$ is
\begin{align*}
  Y_i^2&=0& (1\le i\le n-1),\\
  Y_iY_{i+1}Y_i&=Y_{i+1}Y_iY_{i+1}& (1\le i\le n-2)\\
 Y_jY_i&=Y_iY_j & (j-i\ge 2). 
\end{align*}
In this case, the longest word
is the permutation $(1,n)(2,n-1)\dots$ reversing the $n$ symbols.

The ring $\NW$ has a free $\bZ$-basis consisting of the elements $Y_w$ with $w\in W$,
and multiplication $Y_wY_{w'}=Y_{ww'}$ if $\ell(ww')=\ell(w)+\ell(w')$ and $Y_wY_{w'}=0$ otherwise.
The element $Y_{\mathop{\sf{Id}}}$ is the identity, and the $Y_w$ for $w\ne 1$ are in the nilradical of $\NW$,
since every product of more than $\ell(w_0)$ generators
is equal to zero. 

\begin{remark}\label{rk:NWgraded}
We may regard $\NW$ as a \emph{graded ring}, with the grading given by the length
function. So the generators are in degree one, and $\deg(Y_w)=\ell(w)$.
\end{remark}

Let $J$ be the ideal in $\NW$ generated by the $Y_i$. Then
$J^{\ell(w_0)+1}=0$, so $J$ is the augmentation ideal of $\NW$, it is
nilpotent, and $J=\NW\cap J(\bQ\NW)$ where $J(\bQ\NW)$ is the Jacobson radical.
So we can make $\bZ$ into an $\NW$-module by letting all the $Y_i$
act as zero. The presentation above can then be viewed as the quiver
algebra of a quiver
with one vertex, arrows corresonding to the generators $Y_i$, and
relations as above. 
If $\kk$ is a field, we set $\kNW=\kk\otimes_{\bZ}\NW$, a $\kk$-algebra of dimension
$|W|$, and a unique simple module $\kk$. This is a local basic algebra.

\begin{theorem}
For all $n\ge 1$, $\NW/J^n$ is $\bZ$-free.
The ring $\Ext^*_\NW(\bZ,\bZ)$ is $\bZ$-free, and for any coefficient
field $\kk$ we have 
\[ \Ext^*_{\kNW}(\kk,\kk)\cong \kk\otimes_\bZ\Ext^*_\NW(\bZ,\bZ). \]
\end{theorem}
\begin{proof}
The first statement is the content of
Theorem~11.2 of Benson and Lim~\cite{Benson/Lim:descent}.
The point is that $J^n$ is the $\bZ$-span of the monomials $Y_w$ with
$\ell(w)\ge n$, and this is a $\bZ$-summand of $\NW$.
Localising at any prime ideal in $\bZ$, we then see that Hypothesis~E of that paper
holds. Applying Theorem~2.3 of that paper, we deduce the second statement.
\end{proof}

\begin{remark}\label{rk:psi}
The algebra $\kNW$ is Frobenius, but often not symmetric. We define a
trace map $$\Tr\colon \kNW\to\kk$$ to be the linear map taking $Y_w$ to
zero unless $w=w_0$, in which case $\Tr(Y_{w_0})=1$.
Denote by
$\psi$ the automorphism of $W$ given by conjugation by $w_0$, and use
the same symbol to denote the automorphism of $\kNW$ which takes $Y_w$ to $Y_{w_0ww_0}$
(recall that $w_0^2=1$). In the case of the symmetric group $S_n$, the
automorphism $\psi$ reverses the roles of the numbers $1,\dots,n$,
sending $i$ to $n+1-i$.
\end{remark}

\begin{lemma}
The statements $ww'=w_0$,
$w'\psi(w)=w_0$, and $\psi(w')w=w_0$ are equivalent, and we have
\[ \Tr(ab)=\Tr(b\psi(a))=\Tr(\psi(b)a). \]
\end{lemma}
\begin{proof}
  We have
 $ww'=w_0\Leftrightarrow w_0ww'=1\Leftrightarrow w'w_0w=1\Leftrightarrow
\psi(w')w=w_0w'w_0w=w_0$. Similarly, $ww'w_0=1\Leftrightarrow w'w_0w=1\Leftrightarrow
w'\psi(w)=w'w_0ww_0=w_0$. The displayed equation then follows.
\end{proof}

\begin{proposition}
The map $\psi$ is the Nakayama automorphism of $\kNW$ induced
by the trace map $\Tr\colon\kNW\to\kk$.
\end{proposition}
\begin{proof}
  This follows from the lemma.
\end{proof}

The projective cover $P_\kk$ of $\kk$ is the regular representation,
and is also the injective hull. The socle of the regular representation
is the one dimensional subspace spanned by $Y_{w_0}$.
The Loewy length of the algebra is $\ell(w_0)+1$, and the Loewy and socle
series coincide. In particular, an element $Y_w$ is in the $j$th power
of the radical if and only if $\ell(w)=j$.

The space $\Ext^1_{\kNW}(\kk,\kk)$
has dimension equal to the rank of $W$. So this is $n-1$ in the case of $S_n$.
Also in this case, the representation type is finite for $n=2$,
tame for $n=3$, and wild for $n\ge 4$. If you want to imagine
what the regular representation ``looks like'' in this case, take a
permutohedron and dangle it from a vertex. So for $n=3$ a hexagon,
and for $n=4$ a truncated octahedron, inexpertly rendered in the following diagram.

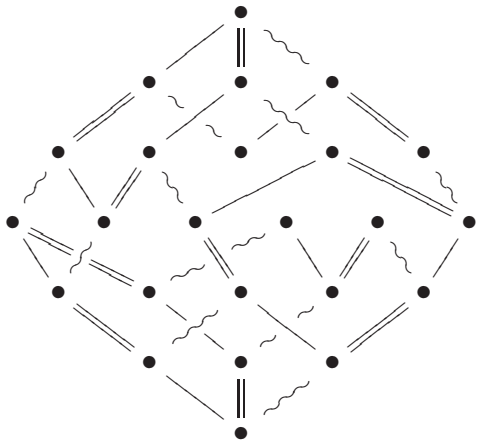
\[ \xymatrix@R=5mm@C=1.5mm{
    &&&&&\bu\ar@{-}[dll]\ar@{=}[d]\ar@{~}[drr]\\
    &&&\bu\ar@{=}[dll]\ar@{~}[drr]|!{[rr];[d]}\hole&&
    \bu\ar@{-}[dll]\ar@{~}[drr]&&
    \bu\ar@{-}[dll]|!{[ll];[d]}\hole\ar@{=}[drr]\\
    &\bu\ar@{-}[dr]\ar@{~}[dl]&&\bu\ar@{=}[dl]\ar@{~}[dr]&&
    \bu\ar@{=}[dr]|!{[rr];[dl]}\hole&&\bu\ar@{-}[dlll]\ar@{=}[drrr]
    &&\bu\ar@{-}[dl]|!{[ll];[dr]}\hole\ar@{~}[dr]\\
    \bu\ar@{-}[dr]\ar@{=}[drrr]|!{[rr];[dr]}\hole&&\bu\ar@{~}[dl]&&\bu\ar@{=}[dr]&&
    \bu\ar@{~}[dlll]|!{[ll];[dl]}\hole\ar@{-}[dr]&&\bu\ar@{=}[dl]\ar@{~}[dr]&&\bu\ar@{-}[dl]\\
    &\bu\ar@{=}[drr]&&\bu\ar@{-}[drr]|!{[rr];[d]}\hole&&\bu\ar@{~}[dll]\ar@{-}[drr]&&
    \bu\ar@{~}[dll]|!{[ll];[d]}\hole&&\bu\ar@{=}[dll]\\
    &&&\bu\ar@{-}[drr]&&\bu\ar@{=}[d]&&\bu\ar@{~}[dll]\\
    &&&&&\bu} \]

See OEIS A008302 for dimensions of Loewy layers, the Mahonian
numbers. The generating function is
\[ \prod_{i=1}^{n-1}(1+x+\dots+x^i) \]
So if $n=4$ then
\[ (1+x)(1+x+x^2)(1+x+x^2+x^3)= 1+3x+5x^2+6x^3+5x^4+3x^5+x^6. \]
The Loewy length is $\binom{n}{2}$.
Here are the first few rows of the triangle of Mahonian numbers
\[ \begin{array}{cccccccccccccccccccccccccccccc}
    1\\
    1&1\\
    1&2&2&1\\
    1&3&5&6&5&3&1\\
    1&4&9&15&20&22&20&15&9&4&1\\
    1&5&14&29&49&71&90&101&101&90&71&49&29&14&5&1\\ 
    1&6&20&49&98&169&259&359&455&531&573&573&531&455&359&259&\cdots\\
    1&7&27&76&174&343&602&961&1415&1940&2493&3017&3450&3736&3836&3736&\cdots\\
     \vdots
\end{array} \]

\section{Canonical form}

Staying with $\NSn$, there is a canonical form for basis elements. 
This section describes this, and gives some technical 
lemmas about it that will be necessary for building the multiple
complex.

It
follows from the factorisation in the group algebra of the symmetric
group:
\[ \sum_{\sigma\in S_n}\sigma = \prod_{i=1}^{n-1}
  (1+s_i+s_{i+1}s_i+s_{i+2}s_{i+1}s_i+\dots+s_{n-1}\dots s_i). \]
Here, the factors are multiplied from left to right as $i$ goes from $1$ to
$n-1$. See for example Theorem~1.1.1 in~\cite{Garsia:2002a}. The
canonical form for $w_0$ is obtained by taking the last term in each factor.
So in
$\NSn$, every basis element has a unique expression as a product
appearing as a term on the right side of the above, with $s_i$
replaced by $Y_i$. For example, $Y_2Y_1Y_4Y_3Y_2$ is in canonical
form, but no word in which $Y_i$ appears more than $i$ times can be in
canonical form. 

\begin{definition}
For $1\le i\le j\le n$ we define 
\[ Y_{[i,j]}=Y_{j-1}Y_{j-2}\dots Y_{i+1}Y_i, \] 
of length $j-i$. In particular, we have $Y_{[i,i]}=1$, $Y_{[i,i+1]}=Y_i$, $Y_{[i,i+2]}=Y_{i+1}Y_i$.
\end{definition}

In terms of this definition, an element in canonical form is as follows:
\[ Y_{[1,m_1]}Y_{[2,m_2]}\dots Y_{[n-1,m_{n-1}]} \]
with each $i\le m_i\le n$.

We can turn the presentation of $\NSn$ into rewriting rules for this
canonical form as follows:
\begin{enumerate}
\item $Y_iY_i\mapsto 0$
\item $Y_jY_i\mapsto Y_iY_j$ ($j-i\ge 2$)
\item $Y_iY_{[i,j]}\mapsto Y_{[i,j]}Y_{i+1}$ ($j-i\ge 2$).
\end{enumerate}

\begin{remark}
The Anick resolution~\cite{Anick:1986a} based on this canonical form is far from
minimal. Our goal is to produce a minimal resolution that is
closely related to the Anick resolution, but using only some of the
terms. This amounts to stripping off an exact free
subcomplex. It is not necessary to understand the Anick resolution to
follow our construction.
\end{remark}

\begin{lemma}\label{le:Yij}
If $j>i$ and $j'>i'\ge i$ then we have
\[ Y_{[i',j']}Y_{[i,j]}=\begin{cases}
Y_{[i,j']}&j=i',\\
Y_{[i,j]}Y_{[i',j']}&j<i',\\
Y_{[i,j]}Y_{[i'+1,j'+1]}&j>j',\\
0&j'\ge j>i'
\end{cases} \]
\end{lemma}
\begin{proof}
This is easy to deduce from the rewriting rules. 
\end{proof}

\begin{definition} 
For $0\le k\le j-i$ we set
\[ Y_{[i,j];k}=Y_{[i,j+1-k]}Y_{[i+1,j+2-k]}\dots Y_{[i+k-1,j]}. \]
This is equal to $Y_w$ where $w$ is the
permutation $(j\,j-1\,\dots\,i)^k$. In particular, we have
$Y_{[i,j];0}=1$, $Y_{[i,j];1}=Y_{[i,j]}$.
\end{definition}

So for example
\begin{align*} 
Y_{[2,7];1} &=Y_6Y_5Y_4Y_3Y_2\\
Y_{[2,7];2}&=(Y_5Y_4Y_3Y_2)(Y_6Y_5Y_4Y_3)\\
Y_{[2,7];3}&=(Y_4Y_3Y_2)(Y_5Y_4Y_3)(Y_6Y_5Y_4)\\
Y_{[2,7];4}  &=(Y_3Y_2)(Y_4Y_3)(Y_5Y_4)(Y_6Y_5) \\
Y_{[2,7];5} &= Y_2Y_3Y_4Y_5Y_6.
\end{align*}
correspond to the powers of the permutation $(7\,6\,5\,4\,3\,2)$.

\begin{lemma}\label{le:Yijk}
If $1\le k\le j-i$ and $j'>i'\ge i$ then
\[ Y_{[i',j']}Y_{[i,j];k}=\begin{cases}
Y_{[i,i'-1];i'+k-j-1}Y_{[i'+k+1-j-i,j']}Y_{[i'+k+2-j-i,j];j-i'}&j\ge i'\ge j+1-k\\
Y_{[i,j];k}Y_{[i',j']}& j<i'\\
Y_{[i,j];k}Y_{[i'+k,j'+k]}&j+1-k>j'\\
0&j'\ge j+1-k>i'
\end{cases} \]
\end{lemma}
\begin{proof}
This follows from Lemma~\ref{le:Yij} by induction on $k$.
\end{proof}

\begin{lemma}\label{le:Y}
We have the following relations among the $Y_{[i,j];k}$.
\begin{enumerate}
\item If $i\le i'<j'$ and $j-j'\ge k$ then we have
\[  Y_{[i',j'];k'}Y_{[i,j];k}=
 Y_{[i,j];k}Y_{[i'+k,j'+k];k'} \]
\item If $0< k' < k \le j-i$ then we have 
\[Y_{[i,j];k}Y_{[i,i+k-1];k'}=Y_{[i,j];k'}Y_{[i+k',j];k-k'}. \]
\end{enumerate}
\end{lemma}
\begin{proof}
Both parts follow from Lemma~\ref{le:Yijk}
 by induction on $k'$.  
\end{proof}

\begin{lemma}\label{le:YYij}
Let $Y$ be a monomial in $\NSn$. 
\begin{enumerate}
\item If $i\le i'<j'<j$ and $YY_{[i,j]}\ne 0$ then $YY_{[i',j']}\ne
0$ and $YY_{[i,j]}Y_{[i',j']}\ne 0$.
\suspend{enumerate}
Suppose that either $i\le i'<j'<j$
or $[i',j']$ is disjoint from $[i,j]$.
\resume{enumerate}
\item If $YY_{[i,j]}\ne 0$ and $YY_{[i',j']}\ne 0$ and
$YY_{[i,j]}Y_{[i',j']}\ne 0$.
\item If  there exists $Y'$ such that $Y'Y_{[i',j']}=Y$ and there exists
  $Y''$ such that $Y''Y_{[i,j]}=Y$ then there exists $Y'''$ such
  that $Y'''Y_{[i',j']}Y_{[i,j]}=Y$.
\end{enumerate}
\end{lemma}
\begin{proof}
Put $Y$ in canonical form, and then repeatedly use Lemma~\ref{le:Yij}.
\end{proof}

\begin{lemma}\label{le:YYijk}
Let $Y$ be a monomial in $\NSn$. 
Suppose that either $i\le i'<j'$ and $j-j'\ge k$, or $[i',j']$ is disjoint from $[i,j]$.
\begin{enumerate}
\item 
If $YY_{[i,j];k}\ne 0$ and $YY_{[i',j'];k'}\ne 0$ then
$YY_{[i',j'];k'}Y_{[i,j];k}\ne 0$.
\item 
If  there exists $Y'$ such that $Y'Y_{[i',j'];k}=Y$ and there exists
  $Y''$ such that $Y''Y_{[i,j];k}=Y$ then there exists $Y'''$ such
  that $Y'''Y_{[i',j'];k'}Y_{[i,j];k}=Y$.
\end{enumerate}
\end{lemma}
\begin{proof}
This follows from Lemma~\ref{le:YYij} and induction.
\end{proof}

\begin{remark}
Lemma~\ref{le:YYijk} will eventually be used in order to prove
exactness of our free resolution of $\bZ$ over $\NSn$.
\end{remark}

\section{The ring $Z$}\label{se:Z}

In this section we describe a ring $Z$ that will turn out to be
the ring
$\Ext^*_{\NSn}(\bZ,\bZ)$. This will be used in
the next section to construct a multiple complex that computes this
cohomology ring.

\begin{definition}\label{def:Z}
Let $Z$ be the non-positively graded $\bZ$-free ring generated by
symbols $z_{i,j}$ for $1\le i,j\le n$ $(i\ne j)$ of degree $|z_{i,j}|=-|j-i|$, subject to the following
relations.
\begin{enumerate}
\item $z_{j,i}=(-1)^{j-i}z_{i,j}$,
\item If $[i',j']\subseteq [i,j]$ then
$z_{i,j}z_{i',j'} =(-1)^{(j-i)(j'-i')}z_{i+j-i',i+j-j'}z_{i,j}$. 
\item If $[i,j]$ is disjoint from $[i',j']$ then
$z_{i,j}z_{i',j'} =(-1)^{(j-i)(j'-i')}z_{i',j'}z_{i,j}$.
\item If $[i,j]$ and $[i',j']$ overlap, but neither is contained in the other, then
  $z_{i,j}z_{i',j'}=0$.
\end{enumerate}
\end{definition}

The way to visualise the second of these relations is that the interval $[i',j']$
gets reflected in the interval $[i,j]$ to give the interval
$[i+j-i',i+j-j']$. Note also that if we put $i=i'$, $j=j'$ in the
second relation, and then use the first relation, the signs all cancel to give no information; the
$z_{i,j}$ do not square to zero.

\begin{remark}\label{rk:Zbigraded}
Since $\NSn$ is a graded ring (see Remark~\ref{rk:NWgraded}),
it will be helpful to regard $Z$ as a doubly graded ring, with
$z_{i,j}$ in degree $(-|j-i|,-\frac{(|j-i|)(|j-i|+1)}{2})$. The first grading
is the homological grading as an element of $\Ext^*_{\NSn}(\bZ,\bZ)$, and the second is
the internal grading coming from regarding $\NSn$ as a graded
ring. The differentials in the resolution we construct will then
preserve internal degree. 
\end{remark}

\begin{remark}\label{rk:z*}
The map sending $z_{i,j}$ to itself
and reversing the order of multiplication is an involutory
anti-automorphism of $Z$, so we have $Z\cong Z^\op$. We write $z^*$
for the image of $z\in Z$ under this antiautomorphism.

The map sending $z_{i,j}$ to $(-1)^{(n-1)(i-j)}z_{n+1-i,n+1-j}$ and preserving the
order of multiplication is an involutory automorphism of $Z$ which we
denote $z\mapsto z^\dagger$. This
satisfies $z_{1,n}z=z^\dagger z_{1,n}$ and $z_{1,n}=z^\dagger_{1,n}$.
\end{remark}

\begin{definition}\label{def:tZ}
We define another ring $\tilde Z$ to be the non-positively graded
$\bZ$-free ring generated by the symbols $\tz_{i,j}$ for $1\le i,j\le n$
$(i\ne j)$ 
of degree $-|j-i|$, subject to the relations of
Definition~\ref{def:Z} but without the signs.
\end{definition}

\begin{definition}\label{def:Zcanonical}
Using these relations, we can put every non-zero monomial in the
$\tz_{i,j}$ in $\tilde Z$ (or monomial in the $z_{i,j}$ in $Z$) uniquely in
\emph{canonical form}, by which we mean that the $j-i$ are positive,
(non-strictly) increasing, each interval $[i,j]$ appearing as a
subscript either contains or is disjoint from each of the previous
ones,  and consecutive ones $\tz_{i,j}\tz_{i',j'}$ of the
same degree are lexicographic, in the sense that they satisfy either
$i=i'$ or $j<i'$.
We define a non-zero monomial
to be in \emph{reversed canonical form} if it has been obtained from
the canonical form by applying the commutation relations of
Definition~\ref{def:Z}\,(2), (3) to reverse
the terms.
\end{definition}

Going from canonical form to reverse canonical form involves a series
of reflections in various intervals. For example, if $n=9$ then moving
terms across one at a time we have 
\begin{align*}
 \tz_{5,6}\tz_{2,4}\tz_{5,7}\tz_{5,7}\tz_{1,4}\tz_{5,8}\tz_{1,9}
&  =\tz_{1,9}\ \tz_{5,4}\tz_{8,6}\tz_{5,3}\tz_{5,3}\tz_{9,6}\tz_{5,2} \\
&=\tz_{1,9}\tz_{5,2}\ \tz_{2,3}\tz_{8,6}\tz_{2,4}\tz_{2,4}\tz_{9,6}\\
  &=\tz_{1,9}\tz_{5,2}\tz_{9,6}\ \tz_{2,3}\tz_{7,9}\tz_{2,4}\tz_{2,4}\\
&=\tz_{1,9}\tz_{5,2}\tz_{9,6}\tz_{2,4}\ \tz_{4,3}\tz_{7,9}\tz_{2,4}\\
  &=\tz_{1,9}\tz_{5,2}\tz_{9,6}\tz_{2,4}\tz_{2,4}\ \tz_{2,3}\tz_{7,9}\\
  &=\tz_{1,9}\tz_{5,2}\tz_{9,6}\tz_{2,4}\tz_{2,4}\tz_{7,9}\tz_{2,3}.
\end{align*}
A similar computation holds in $Z$, but with signs.
                                                                   
\begin{theorem}\label{th:Z}
There is a canonical bijection between the canonical form basis
elements of degree $-d$ in $\tilde Z$ (or in $Z$) 
and the $(n-1)$-tuples of non-positive integers
$(-t_1,\dots,-t_{n-1})$ with $\sum_{i=1}^{n-1} t_i=d$.
\end{theorem}
\begin{proof}
We define a bijection $f$ from the $(n-1)$-tuples to the canonical
form basis elements by induction on negated degree. We begin the induction by
setting $f(0,\dots,0)=1$. For the inductive step,
suppose that we have
already dealt with negated degrees smaller than $d$.
Given an $(n-1)$-tuple $(-t_1,\dots,-t_{n-1})$, write the set $\{i\mid t_i<0\}$
as a union of intervals $[i_1,j_1],\dots,[i_r,j_r]$ with the beginning
of each interval $i_s$ separated from the end of the previous one
$j_{s-1}$ by at least two, so that $i_s-j_{s-1}\ge 2$ for
$2\le s\le r$. Then we look at $f(-t'_1,\dots,-t'_{n-1})\tz_{i_1,j_1+1}\dots\tz_{i_r,j_r+1}$,
where $t'_i=\max(0,t_i-1)$. This is almost in canonical form, except
for the lexicographic condition. Reordering the generators of the same
degree lexicographically, we define the result to be
$f(-t_1,\dots,-t_{n-1})$.
This inductive process is clearly reversible, so it defines a bijection.
\end{proof}

Let us illustrate this process with the $8$-tuple
$(2,3,3,1,5,4,2,1)$ that corresponds to the word we used as an example
above. We have
\begin{align*}
  f(-2,-3,-3,-1,-5,-4,-2,-1)&=f(-1,-2,-2,0,-4,-3,-1,0)\tz_{1,9}\\
&=f(0,-1,-1,0,-3,-2,0,0)\tz_{1,4}\tz_{5,8}\tz_{1,9}\\
&=f(0,0,0,0,-2,-1,0,0)\tz_{2,4}\tz_{5,7}\tz_{1,4}\tz_{5,8}\tz_{1,9}\\
&=f(0,0,0,0,-1,0,0,0)\tz_{5,7}\tz_{2,4}\tz_{5,7}\tz_{1,4}\tz_{5,8}\tz_{1,9}\\
&=f(0,0,0,0,-1,0,0,0)\tz_{2,4}\tz_{5,7}\tz_{5,7}\tz_{1,4}\tz_{5,8}\tz_{1,9}\\
&=\tz_{5,6}\tz_{2,4}\tz_{5,7}\tz_{5,7}\tz_{1,4}\tz_{5,8}\tz_{1,9}.
\end{align*}
The second to last step here illustrates the need for lexicographic reordering.

\begin{corollary}\label{co:PS}
  The ring $Z$ has Poincar\'e series
  \[ \sum_{i=0}^\infty t^i\rank_\bZ(Z_{-i}) = \frac{1}{(1-t)^{n-1}}. \]
\end{corollary}

\begin{remark}\label{rk:comm}
The subring of $Z$ or of $\tilde Z$ generated by the squares of the
generators is commutative (but not central), 
because reflecting twice gets you back to the original.
The whole ring is finitely generated as
a module over this subring. So $Z$ and $\tilde Z$ are Noetherian affine PI
rings. The theory of PI rings is extensively discussed in
Chapter~13 of~\cite{McConnell/Robson:2001a}.

Given a monomial $z$ in canonical form, the element $z^*$ of Remark~\ref{rk:z*} has the
property that $zz^*$ is a monomial whose canonical form is the double
of that of $z$ in the sense that every term is repeated twice.  These
monomials $zz^*$ therefore all commute and are non-nilpotent. 
\end{remark}

We recall the following definition.

\begin{definition}
A ring $R$ is \emph{semiprime} if for all $a\ne 0$ in $R$,
we have $aRa\ne 0$, and $R$ is \emph{prime} if for all $a\ne 0$, $b\ne 0$ in
$R$, we have $aRb\ne 0$. For a graded ring, these conditions need only
be tested on homogeneous elements.
\end{definition}

\begin{theorem}\label{th:semiprime}
The ring $Z$ is semiprime. If $n\ge 4$ then $Z$ is not prime.
\end{theorem}
\begin{proof}
Let $z$ be a non-zero element of $Z$. Lexicographically order the
monomials in canonical form, where the generators are ordered first by
degree and then by first index, and let $z_0$ be the the earliest monomial with non-zero
coefficient in $z$. Then $z_0z_0^*z_0$ is $z_0$ with each term
repeated three times, and has non-zero coefficient in $zz^*z$. It
follows that $zz^*z\ne 0$. Thus $Z$ is semiprime.

If $n\ge 4$, the elements $z_{1,2}z_{n-1,n}$ and $z_{2,n-1}$ are
non-zero, but $(z_{1,2}z_{n-1,n})Zz_{2,n-1}=0$. So $Z$ is not prime.
\end{proof}

Finally, we shall need the following easy to verify technical fact about the
bijection $f$ above. Recall from Remark~\ref{rk:Zbigraded} that the
internal degree of $\tz_{i,j}$ is $-\frac{(|j-i|)(|j-i|+1)}{2}$. 

\begin{proposition}\label{pr:sum-of-degrees}
Let 
\begin{align*}
\tz&=f(-t_1,\dots,-t_{n-1})\\
\tz'&=f(-t_1,\dots,-t_{i-1},-(t_i+1),\dots,-(t_{j-1}+1),-t_j\dots,-t_{n-1})
\end{align*}
 with $1\le i<j\le n$. Then the internal
degree of $\tz'$ is greater
than or equal to the sum of the internal degrees of $\tz_{i,j}$ and
$\tz$ (i.e., the negated degrees are less than or equal), and the following are equivalent.
\begin{enumerate}
\item
$\tz_{i,j}\tz\ne 0$
\item
$\tz_{i,j}\tz=\tz'$
\item
The internal degree of $\tz'$ is the
sum of the internal degrees of $\tz_{i,j}$ and $\tz$.
\item
For all $i\le k < j$ we have
$t_{i-1}\le t_i\le t_k\ge t_{j-1}\ge t_j$.\qed
\end{enumerate}
\end{proposition}

\section{The multiple complex}

We now describe a multiple complex $C$ indexed by $(n-1)$-tuples of
non-negative integers. Each term is a copy of the regular
representation $P$ of $\NSn$, and there is a corresponding monomial
in $Z$ given by negating and applying Theorem~\ref{th:Z}. We may identify $C$
with $\NSn\otimes \check{Z}$, where $\check{Z}$ is the graded dual of $Z$,
a non-negatively graded ring with a free $\bZ$-basis consisting of the duals of
monomials in $Z$. We write the dual of a monomial in the obvious way,
so that for example $\vz_{1,2}\vz_{1,3}$ is the dual of the
monomial $z_{1,2}z_{1,3}$, even though it's not a product.
The element $\vz_{i,j}$ in $\check{Z}$ will
correspond to the term 
\begin{equation}\label{eq:Anick}
\underline{i\ i+1}\underline{\vphantom{i+1}\overline{\ i}}
  \overline{\ i+2\ i+1\ \underline{i}}\ \dots\ 
\underline{\overline{i}\ j-1\ j-2\ \dots\ i+1\ i}
\end{equation}
in the Anick resolution (so for example $\vz_{1,5}$ corresponds to
$\underline{12}\underline{\overline{1}}\overline{32\underline{1}}\underline{4321}$).

To describe $C$, we first construct a multi-graded object $\tilde
C=\NSn\otimes\check Z$ indexed by
$(n-1)$-tuples of non-negative integers as above, with maps $\td_k$ that
decrease the $k$th coordinate by one. But in contrast with a multiple
complex, the maps $\td_k$ square to zero and commute. The maps in
$\tilde C$ involve no signs. To obtain
$C$ from $\tilde C$, we attach signs in the usual way to make the
differentials $d_k$ square to zero and anticommute, in the sense that
for $k\ne k'$ we have $d_kd_{k'}+d_{k'}d_k=0$.

We first define the maps $\td_k$ on the elements $Y\vz_{i,j}$
(meaning $Y\otimes \vz_{i,j}$) for $Y\in \NSn$.

\begin{definition}
Set
\[ \td_k(Y\vz_{i,j})=\begin{cases}
YY_{[i,j];1+k-i}\,\vz_{i,k}\vz_{1+k,j}&i\le
    k<j\\0& k<i\mbox{ or }k\ge j.\end{cases} \]
Here, $\vz_{i,i}$ and $\vz_{j,j}$ are interpreted as $\check{1}$ in degree zero.
\end{definition}

\begin{example}
We have $\td_i(\vz_{i,i+1})=Y_i$, 
$\td_1(\vz_{1,3})=Y_2Y_1\vz_{2,3}$,
$\td_2(\vz_{1,3})=Y_1Y_2\vz_{1,2}$. Notice that
$\td_1\td_2(\vz_{1,3})=\td_1(Y_1Y_2\vz_{1,2})=Y_1Y_2Y_1$ while
$\td_2\td_1(\vz_{1,3})=\td_2(Y_2Y_1\vz_{2,3})=Y_2Y_1Y_2$. By the braid
relation, $\td_1$ and $\td_2$ commute.

  For one final example
  \begin{align*}
\td_1(\vz_{1,5})&=Y_4Y_3Y_2Y_1\,\vz_{2,5},\\
\td_2(\vz_{1,5})&=Y_3Y_2Y_1Y_4Y_3Y_2\,\vz_{1,2}\vz_{3,5},\\
\td_3(\vz_{1,5})&=Y_2Y_1Y_3Y_2Y_4Y_3\,\vz_{1,3}\vz_{4,5},\\
\td_4(\vz_{1,5})&=Y_1Y_2Y_3Y_4\,\vz_{1,4}.
\end{align*}
\end{example}

It is easy to extend the definition of $\td_k$ to monomials
that have disjoint intervals $[i,j]$, by applying it only to the
monomial having an
interval with $i\le k<j$, if there is one, and zero otherwise.
But we defer the definition on more general monomials
until later in this section.

\begin{proposition}\label{pr:dz}
The differentials $\td_k$ commute and square to zero on $\vz_{i,j}$.
\end{proposition}
\begin{proof}
If $i\le k<j$ then
\begin{align*}
\td_k\td_k\vz_{i,j}&=Y_{[i,j];1+k-i}\,d_k\vz_{i,k}\vz_{1+k,j}
  =0
\end{align*}               
because $d_k\vz_{i,k}=d_k\vz_{1+k,j}=0$.
  
If $i\le k<k'<j$ then
\begin{align*}
  d_kd_{k'}\vz_{i,j}&=
Y_{[i,j];1+k'-i}
  \,d_k\vz_{i,k'}\vz_{1+k',j}\\
  &  = Y_{[i,j];1+k'-i}Y_{[i,k'];1+k-i}\,
    \vz_{i,k}\vz_{1+k,k'}\vz_{1+k',j} \\
  &=Y_{[i,j];1+k'-i}Y_{[i,k'];1+k-i}\,
    \vz_{i,k}\vz_{1+k,k'}\vz_{1+k',j} 
\end{align*}
 while
\begin{align*}
  d_{k'}d_k\vz_{1,j}&=
Y_{[i,j];1+k-i}
  \,d_{k'}\vz_{i,k}\vz_{1+k,j}\\
  & =Y_{[i,j];1+k-i}Y_{[k+1,j];k'-k}
    \,\vz_{i,k}\vz_{k+1,k'}\vz_{1+k',j}\\
  & =Y_{[i,j];1+k-i}Y_{[k+1,j];k'-k}
    \,\vz_{i,k}\vz_{k+1,k'}\vz_{1+k',j}.
\end{align*}
These are equal because by Lemma~\ref{le:Y}\,(2),
\begin{equation*} 
Y_{[i,j];1+k'-i}Y_{[i,k'];1+k-i}=Y_{[i,j];1+k-i}Y_{[k+1,j];k'-k}. 
\qedhere
\end{equation*}
\end{proof}

Next, we define the action of the differential $\td_k$ on
monomials in $\check{Z}$. A monomial $\vz$ corresponds
to an $(n-1)$-tuple $(t_1,\dots,t_{n-1})$ via negating and using Theorem~\ref{th:Z},
and $\td_k$ takes it to an element of $\NSn$ times the
monomial in $\check{Z}$ obtained by decreasing $t_k$ to $t_k-1$ if $t_k>0$,
and takes it to zero if $t_k=0$.

To determine the element of $\NSn$ multiplying this monomial,
first write $\vz$ in canonical order, and look for the first term
$\vz_{i,j}$ with $i\le k<j$. Then look at the corresponding term
in the reversed canonical order, and see what the corresponding
differential does to that term. That involves an element of $\NSn$,
and this element 
gives the multiple of the appropriate monomial in $Z$.

\begin{remark}
Since $\td_k$ has to do with the involution swapping $k$ and $(k+1)$,
we must take care when reflecting it. If $i\le k<j$ then the
reflection of $\td_k$ in the interval $[i,j]$ is $\td_{i+j-k-1}$ and not $\td_{i+j-k}$.
\end{remark}

\begin{example}
Let $n=3$ and consider the term $\vz_{1,2}\vz_{1,3}$ corresponding to
the tuple $(2,1)$ via Theorem~\ref{th:Z}. Then in reverse canonical
form, this is $\vz_{1,3}\vz_{2,3}$, so to apply $\td_1$, we must
reflect $\td_1$ in the interval $[1,3]$ to get $\td_2$, and then apply
$\td_2$ to $\vz_{2,3}$ to get $Y_2$. The monomial in $\check{Z}$ corresponding
to $(1,1)$ is $\vz_{1,3}$, so we get
\[ \td_1(\vz_{1,2}\vz_{1,3})= Y_2\,\vz_{1,3}. \]
Similarly, for $\td_2$ the relevant term is $\vz_{1,3}$. Since
$\td_2(\vz_{1,3})=Y_1Y_2\,\vz_{1,2}$ and the basis element of $\check{Z}$
corresponding to $(2,0)$ is $\vz_{1,2}^2$, we have
\[  \td_2(\vz_{1,2}\vz_{1,3})=Y_1Y_2\,\vz_{1,2}^2. \]
The reader may wish to check the example
\[ \td_3(\vz_{2,5}\vz_{1,9})=Y_6Y_5Y_7Y_6\vz_{2,3}\vz_{4,5}\vz_{1,9}. \]
\end{example}

\begin{proposition}\label{pr:dzz}
If $[i',j']\subseteq [i,j]$ then the maps $\td_k$ commute
and square to zero on $\vz_{i',j'}\vz_{i,j}$.
\end{proposition}
\begin{proof}
Let us first prove that $d_kd_k(\vz_{i',j'}\vz_{i,j})=0$. Unless
$i'\le k<j'$, the proof is essentially the same as in
Proposition~\ref{pr:dz}. On the other hand, if $i'\le k<j'$ then
\begin{align*}
d_kd_k(\vz_{i',j'}\vz_{i,j})&=d_k(Y_{[i+j-j',i+j-i'];j'-k}\vz_{i',k}\vz_{k+1,j'}\vz_{i,j})\\
&=Y_{[i+j-j',i+j-i'];j'-k}Y_{[i,j];1+k-i}\vz_{i',k}\vz_{k+1,j'}\vz_{i,k}\vz_{k+1,j}.
\end{align*}
The last term of the first is $Y_{[i+j-k-1,i+j-i']}$ while the first term of the
second is $Y_{[i,i+j-k]}$. By Lemma~\ref{le:Yij}, the product of these is zero.

It remains to prove that the $d_k$ commute. We have
\[ \vz_{i',j'}\vz_{i,j}=\vz_{i,j}\vz_{i+j-i',i+j-j'}.\]
If neither $k$ nor $k'$ lies between $i'$ and $j'-1$ then the
computation is essentially the same as in Proposition~\ref{pr:dz},
with the term $\vz_{i',j'}$ carried along.
So let us suppose that $i'\le k'< j'$. The reflection of $\td_{k'}$ in $[i,j]$ is $\td_{i+j-k'-1}$.
Then there are three cases for $k$. If $i\le k<i'$ then
\begin{align*}
\td_{k}\td_{k'}(\vz_{i',j'}\vz_{i,j})
&=Y_{[i+j-j',i+j-i'];i+j-k'-1}\td_{k}(\vz_{i',k'}\vz_{k'+1,j'}\vz_{i,j})\\
&=Y_{[i+j-j',i+j-i'];j'-k'-1}Y_{[i,j];1+k-i}\,\vz_{i',k'}\vz_{k'+1,j'}\vz_{i,k}\vz_{k+1,j}
\end{align*}
while
\begin{align*} 
\td_{k'}\td_{k}(\vz_{i',j'}\vz_{i,j})&=
Y_{[i,j];1+k-i}\td_{k'}(\vz_{i',j'}\vz_{i,k}\vz_{k+1,j})\\ 
&=Y_{[i,j];1+k-i}Y_{[k+1+j-j',k+1+j-i'];j'-k'-1}\,\vz_{i',k'}\vz_{k'+1,j'}\vz_{i,k}\vz_{k+1,j}.
\end{align*}
These expressions are equal by Lemma~\ref{le:Y}\,(1). The other two
cases, namely $i'\le k<j'$ and $j'\le k<j$, are similar.
\end{proof}

\begin{theorem}
The operations $\td_k$ on $\tC$ commute and square to zero.
\end{theorem}
\begin{proof}
The argument is essentially the same as the proofs of
Propositions~\ref{pr:dz} and~\ref{pr:dzz}, but repeated. The details
are tedious, but the essential point is that the relations in
Definition~\ref{def:tZ} for the $\vz_{i,j}$ 
correspond to the relations in Lemma~\ref{le:Y} for the $Y_{[i,j];k}$.
\end{proof}

\begin{lemma}\label{le:Koszul}
If $Y\vz\ne 1$ is a monomial then
\begin{enumerate}
\item Either there exists $i$ with $\td_i(Y\vz)\ne 0$ or there exists $i$ and a
monomial $Y'\vz'$ with $Yz=\td_i(Y'\vz')$.
\item If for $i\ne j$ we have $\td_i(Y\vz)\ne 0$ and $\td_j(Y\vz)\ne 0$ then
we also have $\td_i\td_j(Y\vz)\ne 0$.
\item If for $i\ne j$ we have $Y\vz=\td_i(Y'\vz')$ and $Y\vz=\td_j(Y''\vz'')$ then
  there exists $Y'''\vz'''$ with $\td_j(Y'''\vz''')=Y'\vz'$ and $\td_i(Y'''\vz''')=Y''\vz''$.
\end{enumerate}
\end{lemma}
\begin{proof}
This follows from the definition of the operations
$d_i$ on monomials, using Lemma~\ref{le:YYijk}.
\end{proof}

\begin{definition}
We define the multiple complex $C$ by using the same graded object as
for $\tC$, but with differentials defined as follows. If $\vz\in\check{Z}$ is a
monomial of degree $(t_1,\dots,t_{n-1})$, then
we set
\[ d_k(\vz)=(-1)^{t_1+\dots+t_{k-1}}\td_k(\vz). \]
Since the $\td_k$ square to zero and commute, it follows that the
$d_k$ square to zero and anticommute.
\end{definition}

\begin{remark}\label{rk:cubes}
Lemma~\ref{le:Koszul} may be interpreted as saying that the monomials
$Yz\ne 1$ form a disjoint union of cubes of various dimensions under the
differentials $\td_k$. When we 
change the signs to form the $d_k$, these cubes become exact subcomplexes.
\end{remark}

\begin{remark}\label{rk:Zacts}
The differentials $d_k$ are morphisms of left $\NSn$-modules
since they are given by right multiplications. Given a generator
$z_{i,j}$ in $Z$, the map sending $\vz_{i,j}$ to $\check{1}$ extends to
a map of multiple complexes that permutes the coordinates by
reversing the indices from $i$ to $j$ and involves some signs.
Examining the internal grading of elements (see
Remark~\ref{rk:Zbigraded} and Proposition~\ref{pr:sum-of-degrees}),
this map sends $Y\vz_{i,j}\vz$ to $Y\vz$, 
and sends elements that cannot be written in this form into the
augmentation ideal of the appropriate term (or to zero).
So these maps satisfy the relations in $Z$ if read
modulo the augmentation ideal of $\NSn$, but not strictly. For
example, we have
\[ \xymatrix{Y_1Y_2\vz_{1,2}^3\ar[r]^{z_{1,2}}&Y_1Y_2\vz_{1,2}^2
\ar[r]^{z_{1,2}}&Y_1Y_2\vz_{1,2}\ar[r]^{z_{1,2}}
&Y_1Y_2\vz\ar[r]^(.6){z_{1,2}}&0\\
\vz_{1,2}^2\vz_{1,3}\ar[r]^{z_{1,2}}\ar[u]_{d_2}\ar[d]^{d_1}&
\vz_{1,2}\vz_{1,3}\ar[r]^{z_{1,2}}\ar[u]_{d_2}\ar[d]^{d_1}&
\vz_{1,3}\ar[r]^{z_{1,2}}\ar[u]_{d_2}\ar[d]^{d_1}&
Y_1\vz_{2,3}\ar[r]^(.6){z_{1,2}}\ar[u]_{d_2}\ar[d]^{d_1}&0\\
Y_2\vz_{1,2}\vz_{1,3}\ar[r]^{z_{1,2}}&Y_2\vz_{1,3}\ar[r]^{z_{1,2}}&Y_2Y_1\vz_{2,3}\ar[r]&0} \]
\end{remark}
See Example~\ref{eg:A3} for more examples in a more compact notation.

\section{The cohomology ring}

\begin{theorem}\label{th:exact}
The total complex $\Tot C$, with differential $d=d_1+\dots+d_{n-1}$,
is a free resolution of $\bZ$ as an $\NSn$-module. Tensored with any
field, it produces a minimal resolution.
\end{theorem}
\begin{proof}
The statement that $d^2=0$ follows from the fact that the $d_k$ square
to zero and anticommute. To show that it is a free resolution, we must
show that it is exact in positive degrees. To do this, we construct a
contracting homotopy $h$. Consider a monomial of the form $Yz$.
Let $S(Yz)$ be the set of $i$ with $1\le i < n$ such that
either $d_i(Yz)\ne 0$ or there exists a monomial $Y'z'$ such that
$d_i(\pm Y'z')=Yz$. Lemma~\ref{le:Koszul}\,(1) shows that
$S(Yz)\ne\varnothing$.
Let $m$ be the smallest element of $S(Yz)$. If
$Yz=d_m(\pm Y'z')$ we set $h(Yz)=\pm Y'z'$. Otherwise we set
$h(Yz)=0$. 
Then $h^2=0$, and 
it follows from Lemma~\ref{le:Koszul} and Remark~\ref{rk:cubes}
that for any given 
monomial, either $dh$ or $hd$ acts as the identity, and the other
acts as zero, and so we have $dh+hd=1$.
Finally, minimality follows from the fact that each $d_i$
involves multiplication by an element of the augmentation ideal.
\end{proof}

\begin{corollary}\label{co:Ext}
The ring $Z$ is isomorphic to $\Ext^*_\NSn(\bZ,\bZ)$.
\end{corollary}
\begin{proof}
It follows from Theorem~\ref{th:exact} that $\Ext^*_{\NSn}(\bZ,\bZ)$
is isomorphic to $\Hom_{\NSn}(\Tot C,\bZ)$. This is dual to
$\check{Z}$, hence isomorphic to $Z$.  Multiplications by elements of
$Z$ act as chain maps, as in Remark~\ref{rk:Zacts}, and now that
elements of the augmentation ideal act as zero, the ring structure of $Z$
corresponds to composition of chain maps.
\end{proof}

\section{Examples}\label{se:eg}

In this section, we examine in detail the cases $n=3$ and $n=4$, and
then we describe what happens for some other Coxeter types. In the
following, $Z$ denotes $\Ext_{\NW}(\bZ,\bZ)$. 
If $\kk$ is a field, we write $\kk Z$ for $\kk\otimes_\bZ Z
\cong \Ext^*_{\kNW}(\kk,\kk)$.

\begin{example}\label{eg:A2}
In the case $n=3$ (Dynkin type $A_2$), the 
ring $Z=\Ext^*_{\NSn}(\bZ,\bZ)$ 
has three generators, $x=z_{1,2}$ and $y=z_{2,3}$ of (homological) degree
$-1$ and one generator $z=z_{1,3}$ of degree $-2$. The relations are $xy=yx=0$,
$xz+zy=0$ and $yz+zx=0$. This ring is prime, and has a faithful representation over
$\bZ[t_1,t_2]$ given by
\[ x\mapsto\left(\begin{matrix}t_1&0\\0&0\end{matrix}\right)
\qquad
y\mapsto\left(\begin{matrix}0&0\\0&-t_1\end{matrix}\right)
\qquad
z\mapsto\left(\begin{matrix}0&t_2\\t_2&0\end{matrix}\right). \]
If we tensor with any field and choose non-zero values for $t_1$ and
$t_2$, this gives a surjective map $\kk Z \to \Mat_2(\kk)$, so the
representation is absolutely irreducible.\bigskip
\end{example}

\begin{example}\label{eg:A3}
In the case $n=4$ (Dynkin type $A_3$), the ring $Z=\Ext^*_{\NSn}(\bZ,\bZ)$ has six
generators, $u=z_{1,2}$, $v=z_{2,3}$,
$w=z_{3,4}$, $x=z_{1,3}$, $y=z_{2,4}$, $z=z_{1,4}$ with degrees
$|u|=|v|=|w|=-1$, $|x|=|y|=-2$, $|z|=-3$, and relations
\begin{gather*} 
uv=vu=0,\quad vw=wv=0,\quad uw+wu=0,\\
ux+xv=0,\quad vx+xu=0,\quad vy+yw=0,\quad wy+yv=0, \\
uy=yu=0,\quad wx=xw=0,\\
uz=zw,\quad vz=zv,\quad wz=zu,\\
xy=yx=0,\quad
xz=zy,\quad yz=zx.
\end{gather*}
The minimal primes in $Z$ are $(uw)$ and $(v,x,y)$, and the map
\[ Z \to Z/(uw) \oplus Z/(v,x,y) \]
is injective. The ring $Z/(uw)$ has a
faithful four dimensional representation over $\bZ[t_1,t_2,t_3]$ as follows:
\begin{gather*}
u \mapsto \left(\begin{matrix}
t_1&0&0&0\\
0&0&0&0\\
0&0&0&0\\
0&0&0&0
\end{matrix}\right)\qquad
v\mapsto\left(\begin{matrix}
0&0&0&0\\
0&\!\!-t_1\!\!&0&0\\
0&0&0&0\\
0&0&0&-t_1
\end{matrix}\right)\qquad
w\mapsto\left(\begin{matrix}
0&0&0&0\\
0&0&0&0\\
0&0&t_1&0\\
0&0&0&0
\end{matrix}\right)
\\
x\mapsto\left(\begin{matrix}
0&t_2&0&0\\
t_2&0&0&0\\
0&0&0&0\\
0&0&0&0
\end{matrix}\right)\qquad
y\mapsto\left(\begin{matrix}
0&0&0&0\\
0&0&0&0\\
0&0&0&t_2\\
0&0&t_2&0
\end{matrix}\right)\qquad
z\mapsto\left(\begin{matrix}
0&0&t_3&0\\
0&0&0&t_3\\
t_3&0&0&0\\
0&t_3&0&0
\end{matrix}\right)
\end{gather*}
If we tensor with any field and choose non-zero values for $t_1$,
$t_2$, and $t_3$, this gives a surjective map $\kk Z \to \Mat_4(\kk)$, so the
representation is absolutely irreducible.
The ring $Z/(v,x,y)$ is
generated by $u$, $w$ and $z$, with relations $uw+wu=0$, $uz+zw=0$,
$wz+zu=0$. This has a faithful four dimensional representation over
$\bZ[t_1,t_2,t_3]$ as follows:
\[ u\mapsto \left(\begin{matrix}
t_1&0&0&0\\
0&-t_1&0&0\\
0&0&0&-t_2\\
0&0&-t_2&0
\end{matrix}\right)\qquad
w\mapsto\left(\begin{matrix}
0&t_2&0&0\\
t_2&0&0&0\\
0&0&-t_1&0\\
0&0&0&t_1
\end{matrix}\right)\qquad
z\mapsto \left(\begin{matrix}
0&0&t_3&0\\
0&0&0&t_3\\
t_3&0&0&0\\
0&t_3&0&0
\end{matrix}\right) \]

Continuing with the case $n=4$, here are the first few layers of the complex
$\tC$. A string of digits $i$ as a label on an arrow indicates a product
of elements $Y_i$. So for example $321$ denotes
$Y_3Y_2Y_1$. These represent right module homomorphisms, so for
example in the bottom left of Layer 0, the square commutes because $Y_2Y_1Y_2=Y_1Y_2Y_1$.

\[  \text{\rm Layer 0:}\qquad\qquad
\vcenter{\xymatrix@=7mm{P\ar[dd]^{2}&&P\ar[dd]^{1}\ar[ll]^{21}&&
P\ar[dd]^{2}\ar[ll]^{12}&&P\ar[dd]^{1}\ar[ll]^{21}&&P\ar[dd]^{21}\ar[ll]^{12}\\ \\
P\ar[dd]^{2}&&P\ar[dd]^{1}\ar[ll]^{21}&&P\ar[dd]^{2}\ar[ll]^{12}&&
P\ar[dd]^{12}\ar[ll]^{21}&&P\ar[dd]^{12}\ar[ll]^{2}\\ \\
P\ar[dd]^{2}&&P\ar[dd]^{1}\ar[ll]^{21}&&P\ar[dd]^{21}\ar[ll]^{12}&&
P\ar[dd]^{21}\ar[ll]^{1}&&P\ar[dd]^{21}\ar[ll]^{1}\\ \\
P\ar[dd]^{2}&&P\ar[dd]^{12}\ar[ll]^{21}&&P\ar[dd]^{12}\ar[ll]^{2}&&
P\ar[dd]^{12}\ar[ll]^{2}&&P\ar[dd]^{12}\ar[ll]^{2}\\ \\
P&&P\ar[ll]^{1}&&P\ar[ll]^{1}&&P\ar[ll]^{1}&&P\ar[ll]^{1}}} \]

\[ \text{\rm Layer 1:}\qquad\vcenter{\xymatrix@=7mm{
&P\ar[dd]^{3}\ar[dl]^(.6){23}&&P\ar[dd]^{2}\ar[ll]^{321}\ar[dl]^(.6){123}&&
P\ar[dd]^{3}\ar[ll]^{23}\ar[dl]^(.6){123}&&P\ar[dd]^{2}\ar[ll]^{32}\ar[dl]^(.6){123}&&
P\ar[dd]^{32}\ar[ll]^{23}\ar[dl]^(.6){123}  \\
&&&&&&&&&\\
&P\ar[dd]^{3}\ar[dl]^(.6){23}&&P\ar[dd]^{2}\ar[ll]^{321}\ar[dl]^(.6){123}&&
P\ar[dd]^{3}\ar[ll]^{23}\ar[dl]^(.6){123}&&P\ar[dd]^{23}\ar[ll]^{32}\ar[dl]^(.6){123}&&
P\ar[dd]^{23}\ar[ll]^{3}\ar[dl]^(.6){123}\\
&&&&&&&&&\\
&P\ar[dd]^{3}\ar[dl]^(.6){23}&&P\ar[dd]^{2}\ar[ll]^(.6){321}\ar[dl]^(.6){123}&&
P\ar[dd]^{32}\ar[ll]^{23}\ar[dl]^(.6){123}&&P\ar[dd]^{32}\ar[ll]^{2}\ar[dl]^(.6){123}&&
P\ar[dd]^{32}\ar[ll]^{2}\ar[dl]^(.6){123}\\
&&&&&&&&&\\
&P\ar[dd]^{32}\ar[dl]^(.6){23}&&P\ar[dd]^{2132}\ar[ll]^{321}\ar[dl]^(.6){123}&&
P\ar[dd]^{2132}\ar[ll]^{3}\ar[dl]^(.6){123}&&P\ar[dd]^{2132}\ar[ll]^{3}\ar[dl]^(.6){123}&&
P\ar[dd]^{2132}\ar[ll]^{3}\ar[dl]^(.6){123}\\
&&&&&&&&&\\
&P\ar[dl]^{3}&&P\ar[ll]^{1}\ar[dl]^{3}&&P\ar[ll]^{1}\ar[dl]^{3}&&
P\ar[ll]^{1}\ar[dl]^{3}&&P\ar[ll]^{1}\ar[dl]^{3}\\
&&&&&&&&&}} \]

\[ \text{\rm Layer 2:}\qquad \vcenter{\xymatrix@=7mm{
&P\ar[dd]^{2}\ar[dl]^(.6){32}&&P\ar[dd]^{1}\ar[ll]^{321}\ar[dl]^(.6){21}&&
P\ar[dd]^{2}\ar[ll]^{123}\ar[dl]^(.6){321}&&P\ar[dd]^{1}\ar[ll]^{21}\ar[dl]^(.6){321}&&
P\ar[dd]^{21}\ar[ll]^{12}\ar[dl]^(.6){321}\\
&&&&&&&&&\\
&P\ar[dd]^{2}\ar[dl]^(.6){32}&&P\ar[dd]^{1}\ar[ll]^{321}\ar[dl]^(.6){21}&&
P\ar[dd]^{2}\ar[ll]^{123}\ar[dl]^(.6){321}&&P\ar[dd]^{12}\ar[ll]^{21}\ar[dl]^(.6){321}&&
P\ar[dd]^{12}\ar[ll]^{2}\ar[dl]^(.6){321}\\
&&&&&&&&&\\
&P\ar[dd]^{23}\ar[dl]^(.6){32}&&P\ar[dd]^{12}\ar[ll]^{321}\ar[dl]^(.6){21}&&
P\ar[dd]^{2132}\ar[ll]^{123}\ar[dl]^(.6){321}&&P\ar[dd]^{2132}\ar[ll]^{1}\ar[dl]^(.6){321}&&
P\ar[dd]^{2132}\ar[ll]^{1}\ar[dl]^(.6){321}\\
&&&&&&&&&\\
&P\ar[dd]^{32}\ar[dl]^{2}&&P\ar[dd]^{2132}\ar[ll]^{321}\ar[dl]^{1}&&
P\ar[dd]^{2132}\ar[ll]^{3}\ar[dl]^{1}&&P\ar[dd]^{2132}\ar[ll]^{3}\ar[dl]^{1}&&
P\ar[dd]^{2132}\ar[ll]^{3}\ar[dl]^{1}\\
&&&&&&&&&\\
&P\ar[dl]^{3}&&P\ar[ll]^{1}\ar[dl]^{3}&&P\ar[ll]^{1}\ar[dl]^{3}&&P\ar[ll]^{1}\ar[dl]^{3}&&
P\ar[ll]^{1}\ar[dl]^{3}\\
&&&&&&&&&}} \]

\bigskip

\pagebreak[3]
Here are the maps $z_{1,2}$, $z_{2,3}$, $z_{1,3}$:

\[ \begin{array}{cccc}
&\text{level 0}&\text{level 1}&\text{level 2}\\
z_{1,2}:&
\vcenter{\xymatrix@R=3mm@C=4mm{
&&&&\ar[l]_{2}\\
&&&\ar[l]_{1}&\ar[l]\\
&&\ar[l]_{2}&\ar[l]&\ar[l]\\
&\ar[l]_{1}&\ar[l]&\ar[l]&\ar[l]\\
&\ar[l]&\ar[l]&\ar[l]&\ar[l]}} &
\vcenter{\xymatrix@R=3mm@C=4mm{
&&&&\ar[l]_{3}\\
&&&\ar[l]_{2}&\ar[l]\\
&&\ar[l]_{3}&\ar[l]&\ar[l]\\
&\ar[l]_{21}&\ar[l]&\ar[l]&\ar[l]\\
&\ar[l]&\ar[l]&\ar[l]&\ar[l]}} &
\vcenter{\xymatrix@R=3mm@C=4mm{
&&&&\ar[l]_{2}\\
&&&\ar[l]_{1}&\ar[l]\\
&&\ar[l]_{23}&\ar[l]&\ar[l]\\
&\ar[l]_{32}&\ar[l]&\ar[l]&\ar[l]\\
&\ar[l]&\ar[l]&\ar[l]&\ar[l]}}\\ \\
z_{2,3}:&\ \ 
\vcenter{\xymatrix@R=4mm@C=3mm{
\ar[d]&\ar[d]&\ar[d]&\ar[d]&\ar[d]^1\\
\ar[d]&\ar[d]&\ar[d]&\ar[d]^2&\\
\ar[d]&\ar[d]&\ar[d]^1&\\
\ar[d]&\ar[d]^2&\\
&&}} & \ \ 
\vcenter{\xymatrix@R=4mm@C=3mm{
\ar[d]&\ar[d]&\ar[d]&\ar[d]&\ar[d]^{2}\\
\ar[d]&\ar[d]&\ar[d]&\ar[d]^{3}&\\
\ar[d]&\ar[d]&\ar[d]^2&\\
\ar[d]^{2}&\ar[d]^{132}&\\
&&}} & \ \ 
\vcenter{\xymatrix@R=4mm@C=3mm{
\ar[d]&\ar[d]&\ar[d]&\ar[d]&\ar[d]^{1}\\
\ar[d]&\ar[d]&\ar[d]&\ar[d]^{2}&\\
\ar[d]^3&\ar[d]^2&\ar[d]^{132}&\\
\ar@{}[d]&&\\
&&}} \\
z_{1,3}:&\qquad
\vcenter{\xymatrix@=4mm{
&\ar[dl]\ar@{<->}[dddrrr]&\ar[dl]&\ar[dl]&\ar[dl]\\
&\ar[dl]&\ar[dl]&\ar[dl]&\ar[dl]\\
&\ar[dl]&\ar[dl]&\ar[dl]&\ar[dl]\\
&\ar[dl]&\ar[dl]&\ar[dl]&\ar[dl]\\
&&&&&}} &\qquad
\vcenter{\xymatrix@=4mm{
&\ar[dl]^(0.3){\!\!\!1}\ar@{<->}[dddrrr]&\ar[dl]&\ar[dl]&\ar[dl]\\
&\ar[dl]^(0.3){\!\!\!1}&\ar[dl]&\ar[dl]&\ar[dl]\\
&\ar[dl]^(0.3){\!\!\!1}&\ar[dl]&\ar[dl]&\ar[dl]\\
&\ar[dl]^(0.3){\!\!\!\!12}&\ar[dl]^(0.3){\!\!\!\!12}&\ar[dl]^(0.3){\!\!\!\!12}&\ar[dl]^(0.3){\!\!\!\!12}\\
&&&&&}} & 
\vcenter{\xymatrix@=4mm{
&\ar@{<->}[dddrrr]&\ar[dl]^(0.3){\!\!\!3}&\ar[dl]&\ar[dl]\\
&&\ar[dl]^(0.3){\!\!\!3}&\ar[dl]&\ar[dl]\\
&&\ar[dl]^(0.3){\!\!\!\!32}&\ar[dl]^(0.3){\!\!\!\!32}&\ar[dl]^(0.3){\!\!\!\!32}\\
&&&&\\
&&&&}}
\end{array} \]
\end{example}\ 

Next we describe the $\Ext$ algebra for the nilCoxeter algebras of type
$B_2$, $G_2$, $I_2(n)$, $B_3$, $H_3$ and $D_4$ for comparison.

\begin{example}\label{eg:B2}
Consider the Coxeter group of type $B_2$. The corresponding nilCoxeter
algebra has generators $Y_1$ and $Y_2$ and relations
\[ Y_1^2=Y_2^2=0,\quad Y_1Y_2Y_1Y_2=Y_2Y_1Y_2Y_1. \]
Consider also the Coxeter group of type $G_2$. The corresponding
nilCoxeter algebra in this case has generators $Y_1$ and $Y_2$ and relations
\[ Y_1^2=Y_2^2=0,\quad Y_1Y_2Y_1Y_2Y_1Y_2=Y_2Y_1Y_2Y_1Y_2Y_1. \]
The Ext algebras in these two cases are isomorphic. In both cases
there are generators $x=z_{1,2}$, $y=z_{2,3}$ of degree $-1$
and $z=z_{1,3}$ of degree $-2$, with relations $xy=yx=0$, $xz+zx=0$,
$yz+zy=0$. This time there are two minimal primes, namely $(x)$ and
$(y)$. Each of $Z/(x)$ and $Z/(y)$ has a faithful matrix
representation over $\bZ[t_1,t_2]$:
\[ y\text{ (or $x$)}\mapsto
  \left(\begin{matrix}t_1&0\\0&-t_1\end{matrix}\right),
\qquad
x\text{ (or $y$)}
\mapsto\left(\begin{matrix}0&0\\0&0\end{matrix}\right)
\qquad
z\mapsto\left(\begin{matrix}0&t_2\\t_2&0\end{matrix}\right). \]
If we tensor with a field $\kk$ of odd characteristic
and choose non-zero values for $t_1$
and $t_2$ we get surjective maps from $\kk Z/(x)$ and $\kk Z/(y)$ to
$\Mat_2(\kk)$. The prime $2$ is a torsion prime for $B_2$ and $G_2$,
so we might expect different behaviour. And indeed, in characteristic
two the ring $\kk Z$ is commutative, and the rings $\kk Z/(x)$ and
$\kk Z/(y)$ are each just a polynomial ring in two variables.
\end{example}

\begin{example}\label{eg:I2(n)}
The Coxeter group of type $H_2$ has nilCoxeter algebra with generators
$Y_1$ and $Y_2$ and relations
\[ Y_1^2=Y_2^2=0,\quad Y_1Y_2Y_1Y_2Y_1=Y_2Y_1Y_2Y_1Y_2. \]
The $\Ext$ ring in this case is isomorphic to that of type $A_2$
discussed in Example~\ref{eg:A2}. More generally, the nilCoxeter algebra of
type $I_2(n)$ 
has the same presentation but where each side of the braid relation
has length $n$ (recall that $I_2(3)=A_2$, $I_2(4)=B_2$, $I_2(5)=H_2$,
$I_2(6)=G_2$). If $n$ is odd the $\Ext$ ring is isomorphic to that
of $A_2$ while if $n$ is even the $\Ext$ ring is isomorphic to that of
$B_2$. 
Recall that conjugation by the longest element of the Coxeter group $I_2(n)$ is
non-trivial for $n$ odd and trivial for $n$ even, which explains the
difference in behaviour of the $\Ext$ rings in these two cases.
The exact value of $n$ may be recoverd by examining the higher Massey
products. 
The Massey product $\langle x,y,x,\dots\rangle$ ($n$ terms)
is the first non-zero one, and is equal to $\langle
y,x,y,\dots\rangle$ ($n$ terms) and to $z$, up to signs.
\end{example}

\begin{example}\label{eg:B3}
Consider the Coxeter group of type $B_3$. The corresponding
nilCoxeter algebra has generators $Y_1$, $Y_2$, $Y_3$ and relations
\[ Y_1^2=Y_2^2=Y_3^2=0,\ Y_1Y_2Y_1=Y_2Y_1Y_2,\
  Y_2Y_3Y_2Y_3=Y_3Y_2Y_3Y_2,\ Y_1Y_3=Y_3Y_1. \] 
The $\Ext$ algebra has generators $u=z_{1,2}$, $v=z_{2,3}$,
$w=z_{3,4}$ of degree $-1$, $x=z_{1,3}$, $y=z_{2,4}$ 
of degree $-2$, and $z=z_{1,4}$ of degree $-3$, with relations
\begin{gather*}
uv=vu=0,\quad vw=wv=0,\quad uw+wu=0,\\
ux+xv=0,\quad vx+xu=0,\quad vy+yv=0,\quad wy+yw=0,\\
uy=yu=0,\quad wx=xw=0,\\
uz=zu,\quad vz=zv,\quad wz=zw,\\
xy=yx=0,\quad xz+zx=0,\quad yz=zy.\\
\xymatrix{u\ar@{-}[r]^x&v\ar@{=}[r]^y&w}
\end{gather*}
This is not isomorphic to the $\Ext$ algebra for type $A_3$, and again
the behaviour in characteristic two is different from other characteristics.
\end{example}

\begin{example}\label{eg:H3}
Consider the Coxeter group of type $H_3$. The corresponding nilCoxeter
algebra has generators $Y_1$, $Y_2$, $Y_3$ and relations
\[ Y_1^2=Y_2^2=Y_3^2=0,\ Y_1Y_2Y_1=Y_2Y_1Y_2,\
  Y_2Y_3Y_2Y_3Y_2=Y_3Y_2Y_3Y_2Y_3,\ Y_1Y_3=Y_3Y_1. \]
The $\Ext$ algebra has generators $u=z_{1,2}$, $v=z_{2,3}$,
$w=z_{3,4}$ of degree $-1$, $x=z_{1,3}$, $y=z_{2,4}$ of degree $-2$, and
$z=z_{1,4}$ of degree $-3$, with relations
\begin{gather*}
uv=vu=0,\quad vw=wv=0,\quad uw+wu=0,\\
ux+xv=0,\quad vx+xu=0,\quad vy+yw=0,\quad wy+yv=0,\\
uy=yu=0,\quad wx=xw=0,\\
uz=zu,\quad vz=zv,\quad wz=zw,\\
xy=yx=0,\quad xz+zx=0,\quad yz+yz=0.\\
\xymatrix{u\ar@{-}[r]^x&v\ar@{-}[r]_5^y&w}
\end{gather*}
\end{example}

\begin{example}\label{eg:D4}
Consider the Coxeter group of type $D_4$. The corresponding nilCoxeter
algebra has generators $Y_1$, $Y_2$, $Y_3$ and $Y_4$ and relations
\begin{gather*}
 Y_1^2=Y_2^2=Y_3^2=Y_4^2=0,\ Y_1Y_2Y_1=Y_2Y_1Y_2,\
  Y_1Y_3Y_1=Y_3Y_1Y_3,\\ Y_1Y_4Y_1=Y_4Y_1Y_4,\ Y_2Y_3=Y_3Y_2,\
  Y_2Y_4=Y_4Y_2,\ Y_3Y_4=Y_4Y_3. 
\end{gather*}
The $\Ext$ algebra has generators $p$, $q$, $r$, $s$, $t$, $u$, $v$,
$w$, $x$, $y$, $z$ with $|p|=|q|=|r|=|s|=-1$, 
$|t|=|u|=|v|=-2$, $|w|=|x|=|y|=-3$, $|z|=-4$, with relations
\begin{gather*}
pq=qp=0,\ pr=rp=0,\ ps=sp=0,\ qr+rq=0,\ qs+sq=0,\ rs+sr=0,\\
pt+tq=0,\ qt+tp=0,\ pu+ur=0,\ ru+up=0,\ pv+vs=0,\ sv+vp=0, \\
qu=uq=0,\ qv=vq=0,\ rt=tr=0,\ rv=vr=0,\ st=ts=0,\ su=us=0, \\
pw=wp,\ px=xp,\ py=yp,\ qw=wr,\ rw=wq,\ qx=xs,\ sx=xq,\ ry=ys,\ sy=yr,\\
qy=yq=0,\ rx=xr=0,\ sw=ws=0,\ tu=ut=0,\ tv=vt=0,\ uv=vu=0, \\
pz+zp=0,\ qz+zq=0,\ rz+zr=0,\ sz+zs=0,\\ 
tw+wu=0,\ uw+wt=0,\ 
tx+xv=0,\ vx+xt=0,\ uy+yv=0,\ vy+yu=0,\\ 
ty=yt=0,\ ux=xu=0,\ vw=wv=0,\ tz+zt=0,\ uz+zu=0,\ vz+zv=0, \\
wx=xw=0,\ wy=yw=0,\ xy=yx=0,\ wz=zw,\ xz=zx,\ yz=zy.\\
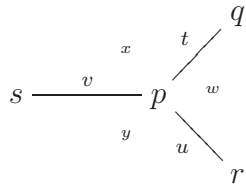
\xymatrix@=6mm{&&&q\ar@{-}[dl]_t\ar@{}[ddl]|{\quad {}^w}\\
s\ar@{-}[rr]^v\ar@{}[urrr]|{{}^x}\ar@{}[drrr]|{{}_y}&&p&\\
&&&r\ar@{-}[ul]^u}
\end{gather*}
Here, $p$, $q$, $r$, $s$ correspond to the vertices, $t$, $u$ and $v$
correspond to the edges, 
$w$, $x$ and $y$ correspond to the illustrated three vertex
subdiagrams, and $z$ corresponds to the whole diagram.
\end{example}

\begin{example}
As a last example, consider the super nilCoxeter algebra $A_3$, as described in Section~8
of Rosso and Savage~\cite{Rosso/Savage:2015a}. This has generators
$Y_1$, $Y_2$, $Y_3$ and relations
\[ Y_1^2=Y_2^2=Y_3^2=0,\quad Y_1Y_2Y_1=Y_2Y_1Y_2,\quad
Y_2Y_3Y_2=Y_3Y_2Y_3,\quad Y_1Y_3+Y_3Y_1=0. \]
The $\Ext$ algebra has generators $u=z_{1,2}$, $v=z_{2,3}$, $w=z_{3,4}$, $x=z_{1,3}$,
$y=z_{2,4}$, $z=z_{1,4}$ with degrees $|u|=|v|=|w|=-1$, $|x|=|y|=-2$,
$|z|=-3$ with relations
\begin{gather*}
uv=vu=0,\quad vw=wv=0,\quad uw=wu,\\
ux+xv=0,\quad vx+xu=0,\quad vy+yw=0,\quad wy+yv=0, \\
uy=yu=0,\quad wx=xw=0, \\
uz+zw=0,\quad vx+xv=0,\quad wz+zu=0, \\
xy=yz=0,\quad xz+zy=0,\quad yz+zx=0.
\end{gather*}
As one might expect, this only differs from Example~\ref{eg:A3} by
some signs. 
\end{example}

\section{PI degree}

In this section, we determine the PI degree of $Z$, defined as
follows. Recall from Posner's theorem~\cite{Posner:1960a} 
(see also 13.3.6 and 13.6.7
of~\cite{McConnell/Robson:2001a}) that a prime PI ring embeds in an
$m\times m$ matrix ring for some $m>0$, over a division ring of
dimension $s^2$ over its centre for some $s>0$, 
and the smallest value of $ms$ is called the PI degree. The ring then
satisfies all the identities satisfied by the matrices of size $ms$
but no smaller.

For a Noetherian semiprime PI
ring there are finitely many minimal primes, 
and we define the PI degree to be the largest PI degree of a
quotient by a minimal prime. Recall from Theorem~\ref{th:semiprime}
that $Z$ is a semiprime PI ring.

\begin{remark}\label{rk:gr-tensor}
For ungraded tensor products of $\bZ$-free rings, 
the PI degrees add, because of the way tensor products of central
simple algebras work. For graded tensor products with the
Koszul sign rule, we need to be more careful, but the PI degree is
bounded by twice the sum, using the following lemma.
\end{remark}

\begin{lemma}\label{le:2x2}
Let $\bZ[t_1,\dots,t_i]$ and $\bZ[t_{i+1},\dots,t_r]$ be commutative
graded polynomial rings with $|t_j|=m_j$. Then the graded tensor
product $\bZ[t_1,\dots,t_i]\otimes\bZ[t_{i+1},\dots,t_r]$ embeds in a
$2\times 2$ matrix ring over a polynomial ring $\bZ[u_1,\dots,u_r]$
with $|u_j|=m_j$.
\end{lemma}
\begin{proof}
We embed the tensor product as follows. If $m_j$ is even then $t_j$ maps to
$\left(\begin{smallmatrix} u_j&0\\0&u_j \end{smallmatrix}\right)$.
If $m_j$ is odd and $1\le j\le i$ then $t_j$ maps to
$\left(\begin{smallmatrix} u_j&0\\0&-u_j \end{smallmatrix}\right)$.
If $m_j$ is odd and $i+1\le j\le r$ then $t_j$ maps to
$\left(\begin{smallmatrix} 0&u_j\\u_j&0 \end{smallmatrix}\right)$.
\end{proof}

\begin{notation}
For $1\le n_1\le n_2\le n$, we have an inclusion of the nilCoxeter subalgebra on the generators
$s_{n_1},s_{n_1+1},\dots,s_{n_2-1}$ to $\NSn$. This induces a
restriction map on $\Ext^*(\bZ,\bZ)$ that sends $z_{i,j}$ to $z_{i,j}$
if $n_1\le i,j\le n_2$ and zero otherwise. This quotient map is
surjective, and so we write
$Z[n_1,n_2]$ for the quotient ring of $Z$ by the ideal generated by the subset of
the generators $z_{i,j}$ not satisfying $n_1\le i,j\le n_2$. In particular, we
have $Z=Z[1,n]$ and $Z[i,i]=\bZ$. We write $Z'=Z'[1,n]$ for the subring of $Z$
generated by all the $z_{i,j}$ with $1\le i<j\le n$, except
$z_{1,n}$. Recall that we have $Z=Z'\otimes \bZ[z_{1,n}]$ with
multiplication given by the rule $z_{1,n}z=z^\dagger z_{1,n}$ for
$z\in Z'[1,n]$, see Remark~\ref{rk:z*}.
\end{notation}

\begin{lemma}\label{le:Z'}
The subring $Z'$ of $Z$ is detected on the quotient rings
$Z[1,i]\otimes Z[i+1,n]$ 
for $1\le i<n$, in the sense that
\[  Z' \to \bigoplus_{1\le i< n} Z[1,i]\otimes Z[i+1,n] \]
is injective.
\end{lemma}
\begin{proof}
The kernel of restriction $Z'\to Z[2,n]=Z[1,1]\otimes Z[2,n]$ is
the ideal generated by the $z_{1,j}$ with $2\le j\le n-1$, while the
kernel of restriction $Z'\to Z[1,n-1]= Z[1,n-1]\otimes Z[n,n]$ is the
ideal generated by the $z_{i,n}$ with  $2\le j\le n-1$. The
intersection of these two kernels is the same as the product, and is
generated by the elements $z_{1,i}z_{j,n}$ with $1<i<j<n$, because if
$n-1\ge i\ge j\ge 2$ then $z_{1,i}z_{j,n}=0$. If $1<i<j<n$, the ideal generated by
$z_{1,i}z_{j,n}$ is detected on the quotients $Z[1,k]\otimes Z[k+1,n]$
with $i\le k<j$.
\end{proof}

\begin{remark}
The analogous statement seems to be true for the other finite Coxeter
groups too. In the case of the $D_4$ Dynkin diagram in 
Example~\ref{eg:D4}, there are three $A_3$ subdiagrams 
and an $(A_1)^3$ subdiagram, and the kernel of restriction to these 
is $(z)$. The $(A_1)^3$ diagram is necessary in order to pick up the
element $qrs$.
\end{remark}

\begin{theorem}\label{th:PI-degree}
For $n\ge 2$, the PI degree of $Z$ is $2^{n-2}$. Over any field $\kk$,
the PI degree of $\kk Z$ is also $2^{n-2}$.
\end{theorem}
\begin{proof}
We prove this by induction on $n$. For $n=2$ the ring is just
$\bZ[z_{1,2}]$, which is commutative, and therefore has PI degree one.
For $n = 3$ and $n=4$ it follows from 
Examples~\ref{eg:A2} and~\ref{eg:A3} that the PI degree is $2^{n-2}$.

Suppose this is true for smaller values of $n$. Then for $1<i<n-1$ the PI
degree of $Z[1,i]\otimes Z[i+1,n]$ is bounded by
$2.2^{i-2}.2^{n-i-2}=2^{n-3}$ (see Remark~\ref{rk:gr-tensor}), while
the PI degrees of $Z[1,n-1]$ and 
of $Z[2,n]$ are exactly $2^{n-3}$.

Consider the restriction maps $\tau\colon Z'\to Z[1,k]\otimes Z[k+1,n]$
and $\tau'\colon Z'\to
Z[1,n-k]\otimes Z[n+1-k,n]$. If we have a map 
\[ \rho\colon Z[1,k]\otimes
Z[k+1,n]\to\Mat_m(\bZ[t_1,\dots,t_{n-2}]) \] 
then we have a map
\[ \rho'\colon Z[1,n-k]\otimes
Z[n+1-k,n]\to\Mat_m(\bZ[t_1,\dots,t_{n-2}]) \] 
given by $\rho'\tau'(z)=\rho\tau(z^\dagger)$ (see Remark~\ref{rk:z*}). 
Since
$z_{1,n}z=z^\dagger z_{1,n}$, it makes sense to 
 combine these to give a map $\hat\rho\colon
Z[1,n]\to\Mat_{2m}(\bZ[t_1,\dots,t_{n-1}])$ by setting
\[ \hat\rho(z_{i,j}) = \left(\begin{matrix} 
    \rho\tau(z_{i,j}) 
& 0 \\ 0 &
\rho'\tau'(z_{i,j}) 
\end{matrix}\right)=
\left(\begin{matrix}
\rho\tau(z_{i,j})&0\\
0&\rho\tau(z_{i,j}^\dagger)
\end{matrix}\right) \]
for $\{i,j\}\ne \{1,n\}$, and 
\[ \hat\rho(z_{1,n})= \left(\begin{matrix} 0 & t_{n-1} \\ t_{n-1} & 0
\end{matrix}\right). \]
We have
$\hat\rho(z_{1,n})\hat\rho(z_{i,j})=\hat\rho(z^\dagger_{i,j})\hat\rho(z_{1,n})$
making $\hat\rho$ a ring homomorphism.

For each matrix representation of $Z[1,k]\otimes Z[k+1,n]$, this way we
obtain a matrix representation of $Z$ of twice the dimension.
Since this subalgebra has PI degree at most
$2^{n-3}$, it follows using Lemma~\ref{le:Z'} and
 induction on $n$, that the PI degree of $Z$ is at
most $2^{n-2}$.  

To show that the PI degree of $Z$ is equal to $2^{n-2}$, over $\bZ$ and after
tensoring with any field, let us suppose
by induction that there is a map 
\[ Z[1,n-1]\to \Mat_{2^{n-3}}(\bZ[t_1,\dots,t_{n-2}]) \]
with the property that after tensoring with any field $\kk$ and
setting $t_i=1$ for all $i$, it gives a
surjective map $\kk Z[1,n-1]\to \Mat_{2^{n-3}}(\kk)$. Let us also
suppose that $z_{1,n-1}$ goes to a non-zero matrix in
$\Mat_{2^{n-3}}(\kk)$.  Then using the
case $k=n-1$ in the construction above, we obtain a map 
\[ Z\to \Mat_{2^{n-2}}(\bZ[t_1,\dots,t_{n-1}]). \] 
After tensoring with $\kk$
and setting $t_i=1$ for all $i$, it gives a map $\kk Z \to
\Mat_{2^{n-2}}(\kk)$. The two maps $\rho\tau$ and $\rho'\tau'$ 
have different kernels, since one sends $z_{1,n-1}$ to a non-zero
matrix and $z_{2,n}$ to zero, and the other is the other way round. 
It follows that they give non-isomorphic absolutely
irreducible representations, so this gives a surjection from $Z'$ onto
the block diagonal matrices with two blocks of size $2^{n-3}$. The
element $z_{1,n}$ goes to the matrix 
$\left(\begin{smallmatrix}0&I\\I&0\end{smallmatrix}\right)$, and this
together with the previous block diagonal matrices generate all
matrices. So by induction, we see that this representation surjects $\kk Z$ onto
$\Mat_{2^{n-2}}(\kk)$ for all fields $\kk$. It follows that the PI
degree over $\bZ$, or over any field cannot be less that $2^{n-2}$.
\end{proof}

\section{The signed nilcactus algebra}

The purpose of this section is to show that the algebra
$Z\cong\Ext^*_{\NSn}(\bZ,\bZ)$ and a signed version of the nilcactus
algebra $X$ are Koszul algebras, Koszul dual to each other. We also
examine a spectral sequence whose 
$E_2$ page is this signed nilcactus algebra, and converging to the nilCoxeter
algebra. 

The cactus group associated to a Coxeter group has been much studied
recently in a wide variety of contexts,  see for example 
\cite{Bellingeri/Chemin/Lebed:2024a,
Bonnafe:2016a,
Chmutov/Glick/Pylyavskyy:2020a,
Etingof/Henriques/Kamnitzer/Rains:2010a,
Genevois,
Halacheva/Kamnitzer/Rybnikov/Weekes:2020a,
Henriques/Kamnitzer:2006a,
Losev:2019a}. This has generators $s_J$ corresponding to the connected
subdiagrams of finite type of the Coxeter diagram, and
relations~\cite[Definition~2.9]{Halacheva/Kamnitzer/Rybnikov/Weekes:2020a}
$s_J^2=1$, $s_Ks_J=s_{\theta_K(J)}s_K$ if $J\subset K$,
$s_Ks_J=s_Js_K$ if $J\cup K$ is disconnected. By analogy with the
nilCoxeter algebra, we define the nilcactus algebra to have a
presentation with respect to generators $X_J$, to which we assign a degree $|J|-1$,
satisfying the same relations 
as the $s_J$, except that the relation $s_J^2=1$ is replaced by
$X_J^2=0$. What we need is a signed version.

In the case of the symmetric group, the signed nilcactus algebra $X$ has
generators $X_{i,j}$ for $1\le i,j\le n$, of degree $|X_{i,j}|=|j-i|-1$,
satisfying the following relations.
\begin{enumerate}
\item $X_{j,i}=(-1)^{j-i-1}X_{i,j}$, $X_{i,i}=0$.
\item $X_{i,j}^2=0$.
\item If $[i',j']\subseteq [i,j]$ then
  $X_{i,j}X_{i',j'}=(-1)^{(j-i-1)(j'-i'-1)}X_{i+j-i',i+j-j'}X_{i,j}$.
\item If $[i,j]$ is disjoint from $[i',j']$ then $X_{i,j}X_{i',j'}=(-1)^{(j-i-1)(j'-i'-1)}X_{i',j'}X_{i,j}$.
\end{enumerate}

In general, quadratic algebras are not necessarily Koszul
algebras. Frequently, to show that particular quadratic
algebras are Koszul, one makes use of Gr\"obner basis techniques. In
our case, the canonical forms of Section~\ref{se:Z} are enough.

\begin{theorem}\label{th:Koszul-dual}
The rings $Z$ and $X$ are Koszul algebras, Koszul dual to
each other.
\end{theorem}
\begin{proof}
The presentation of $Z$ given in Section~4 shows that it is a
quadratic algebra with generators $z_{i,j}$, $i<j$.
The reversed canonical form of Definition~\ref{def:Zcanonical} is a
PBW basis for $Z$ in the sense of Section~5 of
Priddy~\cite{Priddy:1970a}. So it follows from Theorem~5.3 of that
paper, that $Z$ is a Koszul algebra. Taking $X_{i,j}$ to be the dual
basis to $z_{i,j}$ ($i<j$), we see that the perpendicular space of the
quadratic relations in $Z$ is the space of quadratic relations in $X$.
One needs to be careful with the signs here, as the generators are not
in degree one, and we are using graded Koszul duality. Consulting
Theorem~4.6 of~\cite{Priddy:1970a} determines the signs given above.
\end{proof}

There are many ways to put gradings on $X$, but for our context $X$ 
naturally comes with three relevant gradings. There are the homological
grading where each $X_{i,j}$ has degree $-1$, the grading coming from
the homological grading on $Z$, and the grading coming from the
internal grading on $\NSn$, and hence on $Z$ and $X$. Thus the degree
of $X_{i,j}$ is the triple $(-1,|j-i|,\frac{|j-i|(|j-i|+1)}{2})$. The
sum of the first two gives the degree described previously, that makes
$Z$ and $X$ Koszul dual.

{\tiny
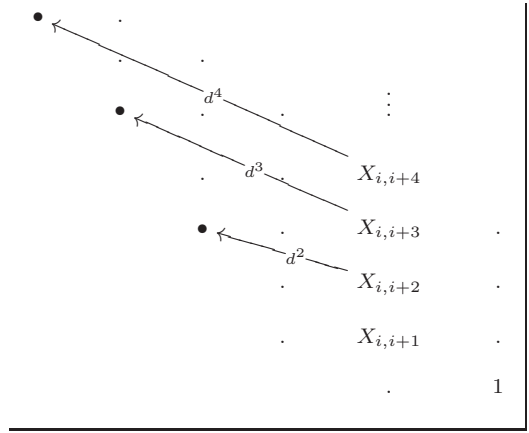
\[ \begin{array}{c|}
\xymatrix@R=2mm@C=7mm{
\bullet&.&&&&\\
&.&.&&&\\
&\bullet&.&.&\vdots&\\
&&.&.&X_{i,i+4}\ar[uuullll]|{d^4}&\\
&&\bullet&.&X_{i,i+3}\ar[uulll]|{d^3}&.\\
&&&.&X_{i,i+2}\ar[ull]|{d^2}&.\\
&&&.&X_{i,i+1}&.\\
&&&&.&1}\\ \\ \hline
\end{array} \]
}
\smallskip

There is an Eilenberg--Moore style spectral sequence $X \Rightarrow
\NSn$ induced by this Koszul duality, where the third degree is
preserved by the differentials. It is illustrated above but with the
third degree suppressed. The $X_{i,i+1}$ are
permanent cycles and correspond to the generators $Y_i$ of $\NSn$. The
element $X_{i,i+2}$ carries a non-zero differential $d^2$ given by 
\[  d^2(X_{i,i+2})=\pm(X_{i,i+1}X_{i+1,i+2}X_{i,i+1}-X_{i+1,i+2}X_{i,i+1}X_{i+1,i+2}) \]
corresponding to the braid relation $Y_iY_{i+1}Y_i=Y_{i+1}Y_iY_{i+1}$.
The element $X_{i,i+3}$ carries a non-zero differential $d^3$ given by
\begin{multline*} 
d^3(X_{i,i+3})=\pm(X_{i,i+2}X_{i+2,i+3}X_{i+1,i+2}X_{i,i+1}+
X_{i+1,i+2}X_{i,i+1}X_{i+1,i+3}X_{i,i+1}\\
\hspace{4em}+X_{i+1,i+2}X_{i+2,i+3}X_{i,i+2}X_{i+2,i+3}
+X_{i+1,i+3}X_{i,i+1}X_{i+1,i+2}X_{i+2,i+3}\\
\hspace{8em}-X_{i+2,i+3}X_{i+1,i+2}X_{i,i+1}X_{i+1,i+3}
-X_{i+2,i+3}X_{i,i+2}X_{i+2,i+3}X_{i+1,i+2}\\
-X_{i,i+1}X_{i+1,i+3}X_{i,i+1}X_{i+1,i+2}
-X_{i,i+1}X_{i+1,i+2}X_{i+2,i+3}X_{i,i+2})
\end{multline*}
expressing the fact that there is a loop in the
rewritings of the socle element $Y_iY_{i+1}Y_iY_{i+2}Y_{i+1}Y_i$ of
$\NS{[i,i+3]}$ using the braid relations, and so on. The first
differential to have non-zero value on $X_{i,j}$ is $d^{|j-i|}$. The
image of this differential represents a sequence of moves going once
around the permutohedron described in Section~\ref{se:Y}.
At the $E_\infty$ page, all that is left are the generators
$X_{i,i+1}$, squaring to zero and now satisfying the braid relations.

\begin{remark}
The cohomology rings of $X$ and of $\NSn$ are both isomorphic to $Z$
as rings. But they are not isomorphic as $A_\infty$-algebras, because
the cohomology of $X$ is formal while that of $\NSn$ is not. So for
example in the cohomology of $\NSn$ we have Massey products
\begin{align*} \langle z_{1,2},z_{2,3},z_{1,2}\rangle &= \pm z_{1,3},\\
 \langle z_{1,3},z_{3,4},z_{2,3},z_{1,2}\rangle&=\pm z_{1,4} 
\\
\langle z_{1,4},z_{4,5},z_{3,4},z_{2,3},z_{1,2}\rangle&=\pm z_{1,5}
\end{align*}
(cf.~\eqref{eq:Anick}) and so on, so that  as an $A_\infty$ algebra it
is generated by degree one elements.
In contrast, in the cohomology of $X$ all Massey products are zero. 

There is an obvious map of algebras $X \to \NSn$, and the induced map
in cohomology is an isomorphism in degree one, sending sending each
$z_{i,i+1}$ to the element with the same name. But this
map in cohomology has to send all generators $z_{i,j}$ with
$|i+j|\ge 2$ to zero, because they can all be expressed as Massey
products of the $z_{i,i+1}$ in the cohomology of $\NSn$. Indeed, there
can be no homomorphism in either direction between $X$ and $\NSn$ that
induces an isomorphism in cohomology, because one is formal and the
other is not.
\end{remark}

\bibliographystyle{amsplain}
\bibliography{../repcoh}

\end{document}